\documentclass[letterpaper, 10 pt, conference]{ieeeconf}  

\IEEEoverridecommandlockouts                              

\usepackage{amsmath}
\usepackage{amssymb}
\usepackage{tikz}
\usepackage{pgfplots}
\usepackage{graphicx}
\usepackage{caption}
\usepackage{subcaption}
\usepackage{xcolor}
\allowdisplaybreaks  
\usepackage{cite}
\usepackage{algorithmicx}
\usepackage{algpseudocode}
\usepackage{algorithm}
\usepackage{subcaption}

\allowdisplaybreaks

\newcommand{\until}[2]{\, U_{[{#1},{#2}]} \, }

\newtheorem{theorem}{Theorem} 
\newtheorem{assumption}{Assumption} 
\newtheorem{problem}{Problem} 
\newtheorem{remark}{Remark}
\newtheorem{lemma}{Lemma}
\newtheorem{corollary}{Corollary}
\newtheorem{example}{Example}

\title{\LARGE \bf Control Barrier Functions for Nonholonomic Systems \\ under  Risk Signal Temporal Logic Specifications }

\author{Lars Lindemann, George J. Pappas, and Dimos V. Dimarogonas
\thanks{This work was supported in part by the Swedish Research Council (VR), the European Research Council (ERC), the Swedish Foundation for Strategic Research (SSF),  the EU H2020 Co4Robots project, the Knut and Alice Wallenberg Foundation (KAW), the DARPA Assured Autonomy program, and the AFOSR grant FA9550-19-1-0265 (Assured Autonomy in Contested Environments).}
\thanks{L. Lindemann and D. V. Dimarogonas are with the Division of Decision and Control Systems, KTH Royal Institute of Technology, 100 44 Stockholm, Sweden. Emails: {\tt\small \{llindem,dimos\}@kth.se}}%
\thanks{G. J. Pappas is with the Department of Electrical and Systems Engineering, University of Pennsylvania, Philadelphia, PA 19104, USA. Email: {\tt\small pappasg@seas.upenn.edu}}%
}

\begin{document}

\maketitle
\thispagestyle{empty}
\pagestyle{empty}

\begin{abstract}
Temporal logics provide a formalism for expressing complex system specifications. A large body of literature has addressed the verification and the control synthesis problem for deterministic systems under such specifications. For stochastic systems or systems operating in unknown environments, however, only the probability of satisfying a specification has been considered so far, neglecting the risk of not satisfying the specification.  Towards addressing this shortcoming, we consider, for the first time, risk metrics, such as (but not limited to) the Conditional Value-at-Risk, and propose risk signal temporal logic. Specifically, we compose risk metrics with stochastic predicates to consider the risk of violating certain spatial specifications. As a particular instance of such stochasticity, we consider control systems in unknown environments and present a determinization of the risk signal temporal logic specification to transform the stochastic control problem into a deterministic one. For unicycle-like dynamics, we then extend our previous work on deterministic time-varying control barrier functions.
\end{abstract}

\section{Introduction}
\label{sec:introduction}

Temporal logic-based control studies the problem of controlling a dynamical system such that a complex specification, expressed as a temporal logic formula, is satisfied. Linear temporal logic (LTL) allows to impose qualitative temporal properties and has been used  in\cite{kress2009temporal,kloetzer2008fully,kantaros2018sampling}. More recently, signal temporal logic (STL) has been considered  \cite{maler2004monitoring}. STL allows to impose quantitative temporal properties, hence being more expressive than LTL. One can additionally associate quantitative semantics with an STL specification which give a real-valued answer to the question whether or not a specification is satisfied, indicating the robustness (severity) of the satisfaction (violation) \cite{donze2,fainekos2009robustness}. Control approaches under STL specifications result in mixed integer linear programs \cite{raman1} or in nonconvex optimization programs \cite{pant2018fly,mehdipour}. Reinforcement learning-based approaches for partially unknown systems have appeared in \cite{varnai2019gradient,aksaray2016q}. An efficient automata-based framework is proposed in \cite{lindemann2019efficient} to decompose the STL specification into STL subspecifications that can be sequentially implemented by low-level feedback control laws, such as those in \cite{lindemann2018control,lindemann2019decentralized} which are based on time-varying control barrier functions. Unlike  other  works, \cite{lindemann2019efficient} directly provides satisfaction guarantees in continuous time. The underlying assumption in \cite{kloetzer2008fully,kantaros2018sampling,raman1,pant2018fly,mehdipour,varnai2019gradient,aksaray2016q,lindemann2019efficient,lindemann2018control,lindemann2019decentralized} is, however,  that the environment is known. For LTL, \cite{yiannis_iros} and \cite{fu2016optimal} assume that the environment is modeled as a semantic map using learning-enabled perception \cite{bowman2017probabilistic} that assign a mean and a variance to each object in the environment. Target beliefs in surveillance games and markov decisions  process-based approaches are respectively presented in \cite{bharadwaj2018synthesis} and \cite{guo2018probabilistic}. Probabilistic computational tree logic and  distribution temporal logic \cite{vasile2016control} account for state distributions and can take chance constraints into account,  but do only consider qualitative temporal properties and do not consider risk metrics as proposed in this work. For STL, literature is sparse and the works in \cite{farahani2018shrinking} and \cite{SadighRSS16} consider chance constraints. 

 Our first contribution is to define risk signal temporal logic (RiSTL) by incorporating risk metrics  \cite{majumdar2020should}, such as (but not limited to)  the Conditional Value-at-Risk \cite{rockafellar2000optimization}, into a temporal logic framework. In particular, we define risk predicates that encode the risk of not satisfying a stochastic STL predicate. On top of these risk predicates, we use the traditional Boolean and temporal operators as in STL.  We also propose quantitative semantics for such specifications.  The second contribution is to show that, under certain conditions, an RiSTL specification can be translated into an STL specification.  We show that these conditions can efficiently be checked for linear predicates, while we argue that, for more general forms, they can be checked numerically. This translation is sound since satisfaction of the STL specification implies satisfaction of the RiSTL specification.  As a particular instance of stochasticity, we consider control systems in unknown environments so that, using this transformation, the stochastic control problem is mapped into a deterministic one. Any existing control method for systems under STL specifications can then be used. We extend, as a third contribution, our previous work on time-varying control barrier functions \cite{lindemann2018control,lindemann2019decentralized} to solve the control problem for unicycle-like dynamics. We emphasize that this is the first work considering unknown environments for continous-time systems under STL alike specifications.


\section{Preliminaries and Problem Formulation}
\label{sec:backgound}

True and false are $\top$ and $\bot$ with $\mathbb{B}:=\{\top,\bot\}$; $\mathcal{N}(\tilde{\boldsymbol{\mu}},\tilde{\Sigma})$ denotes a multivariate normal distribution with mean vector $\tilde{\boldsymbol{\mu}}$ and variance matrix $\tilde{\Sigma}$. Proofs are given in the appendix.

\subsection{Risk Signal Temporal Logic (RiSTL)}
\label{sec:RiSTL}
Let $\boldsymbol{x}:\mathbb{R}_{\ge 0}\to\mathbb{R}^n$ and $\boldsymbol{X}\in\mathbb{R}^{\tilde{n}}$. Signal temporal logic (STL) \cite{maler2004monitoring} is based on signals $\boldsymbol{x}(t)$ and predicates $\mu^{\text{STL}}:\mathbb{R}^n\times \mathbb{R}^{\tilde{n}}\to\mathbb{B}$. Let  $h:\mathbb{R}^n\times \mathbb{R}^{\tilde{n}}\to\mathbb{R}$ be a continuously differentiable function, also called \emph{predicate function}.  A predicate $\mu^{\text{STL}}(\boldsymbol{x}(t),\boldsymbol{X})$ is satisfied at time $t$ if and only if $\boldsymbol{x}(t)$ is such that $h(\boldsymbol{x}(t),\boldsymbol{X})\ge 0$ when $\boldsymbol{X}$ is a deterministic vector. In this paper, however, $\boldsymbol{X}$ is non-deterministic and a random variable. Consider the \emph{probability space} $(\Omega,\mathcal{S}_\Omega,P_\Omega)$  where $\Omega$ is the sample space, $\mathcal{S}_\Omega$ is the Borel $\sigma$-algebra of $\Omega$, and
 $P_\Omega:\mathcal{S}_\Omega\to[0,1]$ is a probability measure. Then $\boldsymbol{X}$ is a measurable function $\boldsymbol{X}:\Omega\to\mathbb{R}^{\tilde{n}}$. Letting $\mathcal{S}$ denote the Borel $\sigma$-algebra of $\mathbb{R}$, the probability space $(\mathbb{R}^{\tilde{n}},\mathcal{S}^{\tilde{n}},P_{\boldsymbol{X}})$ can be associated with $\boldsymbol{X}$ where, for $S\in\mathcal{S}^{\tilde{n}}$, $P_{\boldsymbol{X}}:\mathcal{S}^{\tilde{n}}\to[0,1]$ with $P_{\boldsymbol{X}}(S):=P_\Omega(\boldsymbol{X}^{-1}(S))$ and $\boldsymbol{X}^{-1}(S):=\{\omega\in\Omega|\boldsymbol{X}(w)\in B\}$. Similarly, for a given $\boldsymbol{x}(t)$, one can associate the probability space $(\mathbb{R},\mathcal{S},P)$ with  $h(\boldsymbol{x}(t),\boldsymbol{X})$. We now propose an  extension to STL that  takes  \emph{chance}  and \emph{risk} constraints into account and which we call \emph{risk signal temporal logic} (RiSTL). For a given probability $\delta\in(0,1)$, the truth value of a \emph{chance predicate} $\mu^{\text{Ch}}:\mathbb{R}^n\times \mathbb{R}^{\tilde{n}}\to\mathbb{B}$ at time $t$ is obtained as
\begin{align}\label{eq:prob_predicate}
\mu^{\text{Ch}}(\boldsymbol{x}(t),\boldsymbol{X}):=\begin{cases}
\top & \text{if } P(h(\boldsymbol{x}(t),\boldsymbol{X})\ge 0)\ge \delta\\
\bot &\text{otherwise }
\end{cases}
\end{align}
where $P(h(\boldsymbol{x}(t),\boldsymbol{X})\ge 0)$ denotes the probability that $h(\boldsymbol{x}(t),\boldsymbol{X})\ge 0$, which is the probability of satisfying $\mu^{\text{STL}}(\boldsymbol{x}(t),\boldsymbol{X})$. We further consider \emph{risk predicates} based on risk metrics as advocated in \cite{rockafellar2000optimization,majumdar2020should} and motivated by the fact that chance predicates do not take the left tail of the distribution of $h(\boldsymbol{x}(t),\boldsymbol{X})$ into account. Risk metrics allow to exclude behavior which is deemed more risky than other behavior (see Example \ref{ex:comp} for further motivation). Let  $\mathcal{H}$ denote the set of all random variables derived from $(\Omega,\mathcal{S}_\Omega,P_\Omega)$. Formally, a risk metric is a mapping $R:\mathcal{H}\to\mathbb{R}$.  We are interested in $R(-h(\boldsymbol{x}(t),\boldsymbol{X}))$ to argue about the risk of not satisfying $\mu^{\text{STL}}(\boldsymbol{x}(t),\boldsymbol{X})$. The truth value of a risk predicate $\mu^{\text{Ri}}:\mathbb{R}^n\times \mathbb{R}^{\tilde{n}}\to\mathbb{B}$ at time $t$  is obtained as
\begin{align}\label{eq:risk_predicate}
\mu^{\text{Ri}}(\boldsymbol{x}(t),\boldsymbol{X})&:=\begin{cases}
\top & \text{if } R(-h(\boldsymbol{x}(t),\boldsymbol{X}))\le \gamma\\
\bot &\text{otherwise }
\end{cases}
\end{align} 
for $\gamma\in\mathbb{R}$. Note that $R(\cdot)$ can take different forms with desireable properties such as monotonicity, translational invariance, positive homogeneity, subadditivity, law invariance, or commotone additivity \cite{majumdar2020should}.  The syntax of RiSTL is 
\begin{align}\label{eq:full_RiSTL}
\phi \; ::= \; \top \; | \; \mu^{\text{Ch}}  \; | \; \mu^{\text{Ri}} \; | \; \neg \phi \; | \; \phi' \wedge \phi'' \; | \; \phi'  \until{a}{b} \phi''\;
\end{align}
where $\phi'$ and $\phi''$ are RiSTL formulas and where $\until{a}{b}$ is the until operator with  $a\le b<\infty$. Also define $\phi' \vee \phi'':=\neg(\neg\phi' \wedge \neg\phi'')$ (disjunction), $F_{[a,b]}\phi:=\top \until{a}{b} \phi$ (eventually), and $G_{[a,b]}\phi:=\neg F_{[a,b]}\neg \phi$ (always).  
\begin{remark}
RiSTL allows to impose  specifications like ``the risk of avoiding an obstacle is always less than $\gamma$". It is, however, not possible to impose  ``the risk of always  avoiding an obstacle is less than $\gamma$". While the latter may be more general, we argue that the choice of chance and risk constraints  as in \eqref{eq:prob_predicate} and \eqref{eq:risk_predicate} is more tractable considering that the system in \eqref{system_noise} operates in continuous time and allows to map the stochastic into a deterministic control problem. 
\end{remark}  

Let $(\boldsymbol{x},\boldsymbol{X},t)\models \phi$ denote the satisfaction relation, i.e., if  $\boldsymbol{x}$ satisfies  $\phi$ at $t$ for a particular $\boldsymbol{X}$. We recursively define the RiSTL semantics as $(\boldsymbol{x},\boldsymbol{X},t) \models
 \mu^{\text{Ch}}$ iff	$P(h(\boldsymbol{x}(t),\boldsymbol{X})\ge 0)\ge \delta$, $(\boldsymbol{x},\boldsymbol{X},t) \models
 \mu^{\text{Ri}}$ iff	$R(-h(\boldsymbol{x}(t),\boldsymbol{X}))\le \gamma$, $(\boldsymbol{x},\boldsymbol{X},t) \models \neg\phi$ iff $\neg((\boldsymbol{x},\boldsymbol{X},t) \models \phi)$, $(\boldsymbol{x},\boldsymbol{X},t) \models \phi' \wedge \phi''$ iff $(\boldsymbol{x},\boldsymbol{X},t) \models \phi' \wedge (\boldsymbol{x},\boldsymbol{X},t) \models \phi''$, and $(\boldsymbol{x},\boldsymbol{X},t) \models \phi' \until{a}{b} \phi''$ iff $\exists t'' \in[t+a,t+b] \text{ s.t. }(\boldsymbol{x},\boldsymbol{X},t'')\models \phi'' \wedge \forall t'\in[t,t'']\text{,}(\boldsymbol{x},\boldsymbol{X},t') \models \phi'$.  For a particular $\boldsymbol{X}$, $\phi$ is satisfiable if $\exists \boldsymbol{x}:\mathbb{R}_{\ge 0}\to\mathbb{R}^n$ such that $(\boldsymbol{x},\boldsymbol{X},0)\models \phi$.  Quantitative semantics for RiSTL are denoted by $\rho^{\phi}(\boldsymbol{x},\boldsymbol{X},t)$ and recursively defined as
\begin{align*}
\rho^{\mu^{\text{Ch}}}(\boldsymbol{x},\boldsymbol{X},t)& := P(h(\boldsymbol{x}(t),\boldsymbol{X})\ge 0)-\delta, \\
\rho^{\mu^{\text{Ri}}}(\boldsymbol{x},\boldsymbol{X},t)& :=\gamma-R(-h(\boldsymbol{x}(t),\boldsymbol{X})),\\
 \rho^{\neg\phi}(\boldsymbol{x},\boldsymbol{X},t) &:= 	-\rho^{\phi}(\boldsymbol{x},\boldsymbol{X},t),\\
\rho^{\phi' \wedge \phi''}(\boldsymbol{x},\boldsymbol{X},t) &:= 	\min(\rho^{\phi'}(\boldsymbol{x},\boldsymbol{X},t),\rho^{\phi''}(\boldsymbol{x},\boldsymbol{X},t)),\\
\rho^{\phi' \until{a}{b} \phi''}(\boldsymbol{x},\boldsymbol{X},t) &:= \underset{t''\in [t+a,t+b]}{\max}  \min(\rho^{\phi''}(\boldsymbol{x},\boldsymbol{X},t''), \\ & \hspace{2.5cm}\underset{t'\in[t,t'']}{\min}\rho^{\phi'}(\boldsymbol{x},\boldsymbol{X},t') ), \\
\rho^{G_{[a,b]} \phi}(\boldsymbol{x},\boldsymbol{X},t) &:= \underset{t'\in[t+a,t+b]}{\min}\rho^{\phi}(\boldsymbol{x},\boldsymbol{X},t'),\\
\rho^{F_{[a,b]} \phi}(\boldsymbol{x},\boldsymbol{X},t) &:= \underset{t'\in[t+a,t+b]}{\max}\rho^{\phi}(\boldsymbol{x},\boldsymbol{X},t').\\
\end{align*}
It holds that $(\boldsymbol{x},\boldsymbol{X},t)\models \phi$ if $\rho^\phi(\boldsymbol{x},\boldsymbol{X},t)> 0$ which follows due to \cite[Prop.~16]{fainekos2009robustness}. For $R(\cdot)$, we use, in this paper, the expected value (EV), the Value-at-Risk (VaR), and the Conditional Value-at-Risk (CVaR). The EV of $-h(\boldsymbol{x}(t),\boldsymbol{X})$ is $E[-h(\boldsymbol{x}(t),\boldsymbol{X})]$ which provides a risk neutral risk measure. More risk averse measures are the VaR and the  CVaR as in \cite{rockafellar2000optimization}. The VaR of $-h(\boldsymbol{x}(t),\boldsymbol{X})$ for $\beta\in(0,1)$ is defined as
\begin{align*}
VaR_\beta(-h(\boldsymbol{x}(t),\boldsymbol{X}))&:=\\
&\hspace{-1cm}\min(d\in\mathbb{R}|P(-h(\boldsymbol{x}(t),\boldsymbol{X})\le d)\ge \beta).
\end{align*}
Note in particular that the probability that $-h(\boldsymbol{x}(t),\boldsymbol{X})> VaR_\beta(-h(\boldsymbol{x}(t),\boldsymbol{X}))$ is $1-\beta$. If the cummulative distribution function of $h(\boldsymbol{x}(t),\boldsymbol{X})$ is smooth, as in this case, the CVaR of $-h(\boldsymbol{x}(t),\boldsymbol{X})$ for probability $\beta$ is given by
\begin{align*}
CVaR_\beta(-h(\boldsymbol{x}(t),\boldsymbol{X}))&:=E[-h(\boldsymbol{x}(t),\boldsymbol{X}))|\\
&\hspace{-0.75cm}-h(\boldsymbol{x}(t),\boldsymbol{X}))>VaR_\beta(-h(\boldsymbol{x}(t),\boldsymbol{X}))].
\end{align*}


Risk predicates are fundamentally different from chance predicates and may be advantageous, as illustrated next.

\begin{example}\label{ex:comp}
Let $\boldsymbol{x}:=\begin{bmatrix} x_x & x_y \end{bmatrix}^T\in\mathbb{R}^2$  and $\boldsymbol{X}\sim\mathcal{N}(\tilde{\boldsymbol{\mu}},\tilde{\Sigma})$ with $\boldsymbol{X}:=\begin{bmatrix}\boldsymbol{X}_{O1}^T & \boldsymbol{X}_{O2}^T & \boldsymbol{X}_{R1}^T\end{bmatrix}^T=\begin{bmatrix}{X}_{O1,x} & X_{O1,y} & X_{O2,x} & X_{O2,y} & X_{R1,x} & X_{R1,y}\end{bmatrix}^T\in\mathbb{R}^6$ (see Fig. \ref{fig:exx}).   The uncertainty of $\boldsymbol{X}_{O1}$ and $\boldsymbol{X}_{O2}$ differs in the left and right part of Fig. \ref{fig:exx} and is larger in the right part (see dotted circles). The specification is to always avoid the obstacles indicated by $\boldsymbol{X}_{O1}$ and $\boldsymbol{X}_{O2}$, while eventually reaching the region indicated by $\boldsymbol{X}_{R1}$. Let $\text{UN}\in\{\text{Ch},\text{Ri}\}$ and $\phi^{\text{UN}}:= G_{[0,6]} (\phi_{O1}^{\text{UN}} \wedge \phi_{O2}^{\text{UN}}) \wedge F_{[0,6]} \phi_{R1}^{\text{UN}}$ where 
 \begin{align*}
 \phi_{O1}^{\text{UN}}:=\mu_1^{\text{UN}}\vee \mu_2^{\text{UN}}\vee \mu_3^{\text{UN}}\vee\mu_4^{\text{UN}}\\
 \phi_{O2}^{\text{UN}}:=\mu_5^{\text{UN}}\vee \mu_6^{\text{UN}}\vee \mu_7^{\text{UN}}\vee\mu_8^{\text{UN}}\\
   \phi_{R1}^{\text{UN}}:=\mu_9^{\text{UN}}\wedge \mu_{10}^{\text{UN}}\wedge \mu_{11}^{\text{UN}}\wedge \mu_{12}^{\text{UN}}
 \end{align*}
 encode avoidance of $\boldsymbol{X}_{O1}$ and $\boldsymbol{X}_{O2}$  and reachability of $\boldsymbol{X}_{R1}$, respectively. For $\mu_1^{\text{UN}}$, $\mu_2^{\text{UN}}$, $\mu_3^{\text{UN}}$, and $\mu_4^{\text{UN}}$, define  $h_1(\boldsymbol{x},\boldsymbol{X}):=X_{O1,x}-\epsilon-x_x$, $h_2(\boldsymbol{x},\boldsymbol{X}):=-X_{O1,x}-\epsilon+x_x$, $h_3(\boldsymbol{x},\boldsymbol{X}):=X_{O1,y}-\epsilon-x_y$, $h_4(\boldsymbol{x},\boldsymbol{X}):=-X_{O1,y}-\epsilon+x_y$ for $\epsilon:=0.5$. Define the remaining $h_i(\boldsymbol{x},\boldsymbol{X})$ for  $\mu_i^{\text{UN}}$ with $i\in\{5,\hdots,12\}$ similarly. Now each $\mu_i^{\text{UN}}$  for $i\in\{1,\hdots,12\}$ is either interpreted as a chance predicate ($\text{UN}=\text{Ch})$ with  $\delta_i:=0.5$  or as a risk predicate ($\text{UN}=\text{Ri}$)  with $R(\cdot):=CVaR_{\beta_i}(\cdot)$, $ \gamma_i:=1.5$, and $\beta_i:=0.8$. Fig. \ref{fig:exx} shows the trajectories $\boldsymbol{x}_1(t)$ (blue) and $\boldsymbol{x}_2(t)$ (red). In the left part, both  $\rho^{\phi^{\text{Ch}}}$ and $\rho^{\phi^{\text{Ri}}}$  indicate that $\boldsymbol{x}_2(t)$ satisfies $\phi^\text{UN}$ more (note that $\rho^{\phi^{\text{Ch}}}(\boldsymbol{x}_2,\boldsymbol{X},0)>\rho^{\phi^{\text{Ch}}}(\boldsymbol{x}_1,\boldsymbol{X},0)$ and $\rho^{\phi^{\text{Ri}}}(\boldsymbol{x}_2,\boldsymbol{X},0)>\rho^{\phi^{\text{Ri}}}(\boldsymbol{x}_1,\boldsymbol{X},0)$). The intuition here is that $\boldsymbol{x}_1(t)$ (blue) does not reach the center of $\boldsymbol{X}_{R1}$ as opposed to $\boldsymbol{x}_2(t)$ (red) so that $\boldsymbol{x}_2(t)$ is favoured in both cases since this trajectory satisfies the reachability specification $\phi_{R1}^{\text{UN}}$ better when the uncertainty in $\boldsymbol{X}_{O1}$ and $\boldsymbol{X}_{O2}$ is low. In the right part, however, this uncertainty grows; $\rho^{\phi^{\text{Ch}}}$ still suggests that $\boldsymbol{x}_2(t)$ is the favorable trajectory, while now  $\rho^{\phi^{\text{Ri}}}$, being more risk sensitive, suggest that $\boldsymbol{x}_1(t)$ is more favorable. The reason for this behavior is that the relative importance of the avoidance specifications  $\phi_{O1}^{\text{UN}}$ and $\phi_{O2}^{\text{UN}}$ increases and is more taken into account by the risk predicates. 
\begin{figure}
\centering
\includegraphics[scale=0.2]{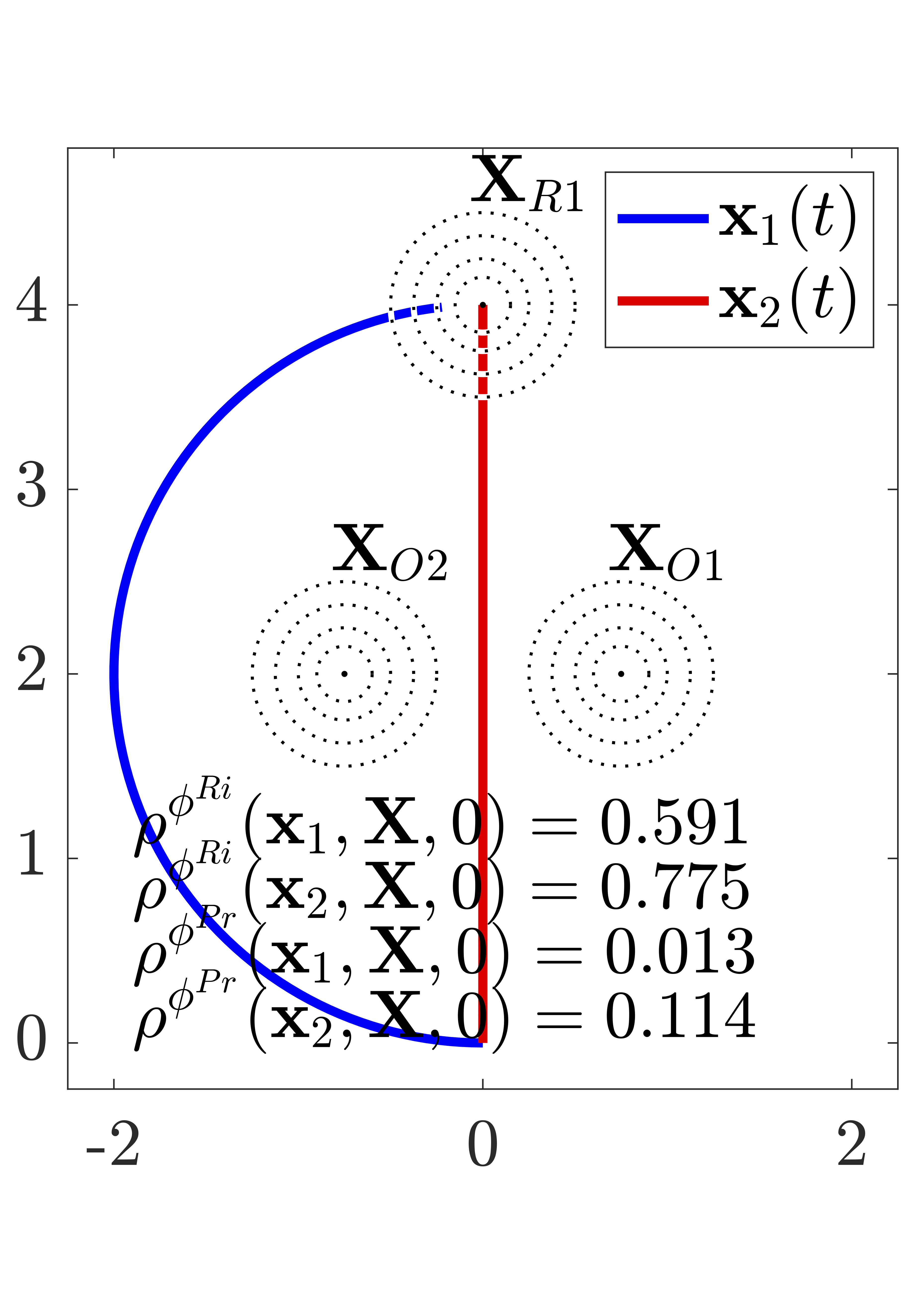}
\includegraphics[scale=0.2]{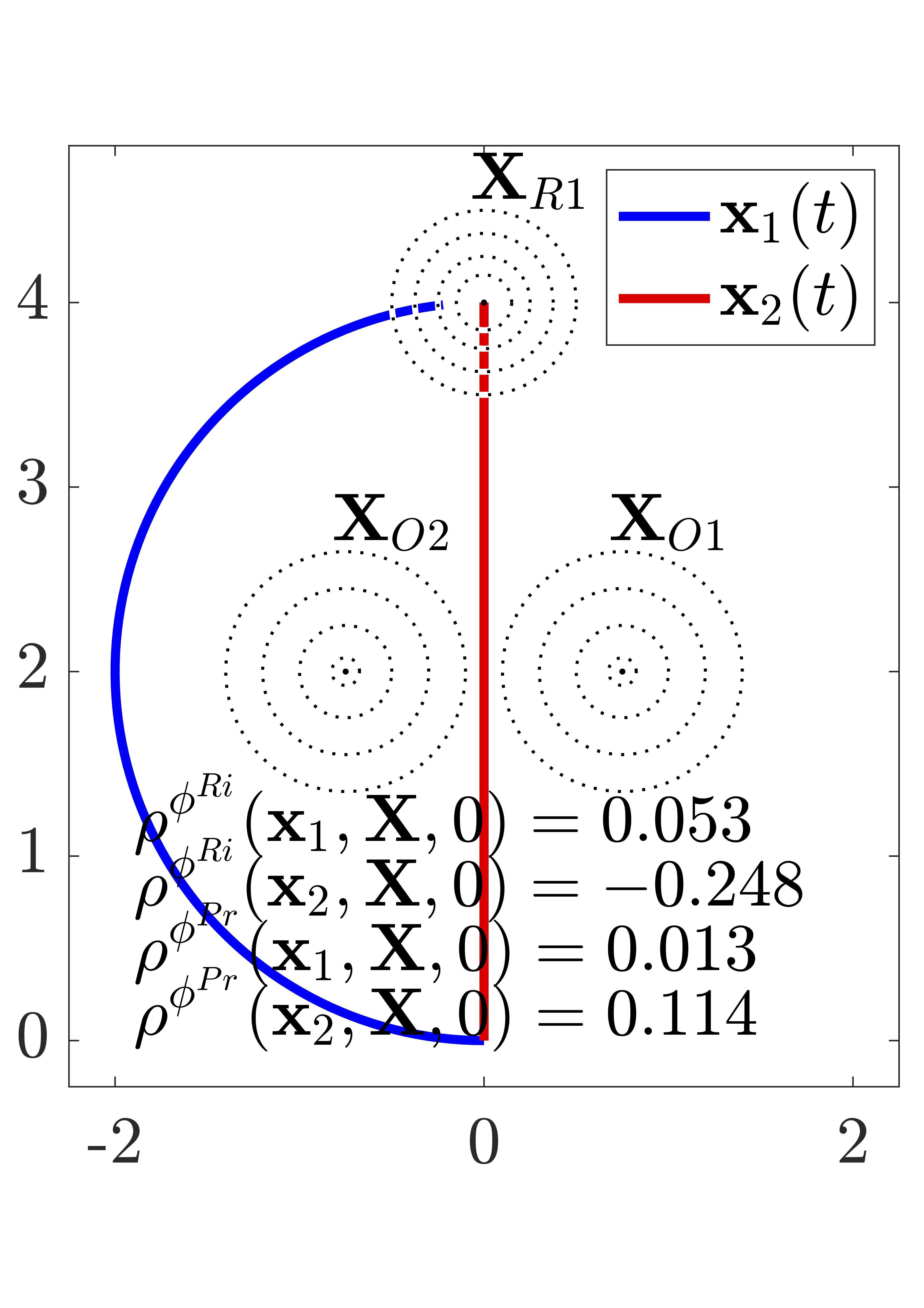}
\caption{Let  $\boldsymbol{X}$ have mean  $\tilde{\boldsymbol{\mu}}:=\big[0.75 \;\, 2 \;\, -0.75 \;\, 2 \;\, 0 \;\, 4\big]^T$ and variance $\tilde{\Sigma}:=0.75\cdot\text{diag}(1,1,1,1,1,1)$ in the left and $\tilde{\Sigma}:=0.75\cdot\text{diag}(2,2,2,2,1,1)$ in the right figure.}
    \label{fig:exx}
\end{figure}
\end{example}


\subsection{Nonholonomic Systems under RiSTL Specifications}
\label{sec:Prob_DEF}

Let $\boldsymbol{z}(t):=\begin{bmatrix} \boldsymbol{x}(t)^T & \theta(t)\end{bmatrix}^T\in\mathbb{R}^3$ where $\boldsymbol{x}(t)$ and $\theta(t)$ are the position and orientation of a unicycle modeled  as  in
\begin{align}\label{system_noise}
\dot{\boldsymbol{z}}(t)=f(\boldsymbol{z}(t))+g(\boldsymbol{z}(t))\boldsymbol{u}+c(\boldsymbol{z}(t),t), \;\boldsymbol{z}(0)\in\mathbb{R}^3 
\end{align}
with control input $\boldsymbol{u}:=\begin{bmatrix} u_1 & u_2 \end{bmatrix}^T\in\mathbb{R}^2$. The functions 
\begin{align*}
f(\boldsymbol{z}):=\begin{bmatrix}
f_{\boldsymbol{x}}(\boldsymbol{z}) \\ f_\theta(\boldsymbol{z})
\end{bmatrix} g(\boldsymbol{z}):=\begin{bmatrix}
\cos(\theta) & 0\\
\sin(\theta) & 0\\
0 & 1
\end{bmatrix} c(\boldsymbol{z},t):=\begin{bmatrix}
c_{\boldsymbol{x}}(\boldsymbol{z},t) \\  c_\theta(\boldsymbol{z},t)
\end{bmatrix}
\end{align*} are  locally Lipschitz continuous in $\boldsymbol{z}$ and piecewise continuous in $t$; $f(\boldsymbol{z})$ is a  known function with bounded $f_\theta(\boldsymbol{z})$ while $c(\boldsymbol{z},t)$ is unknown but bounded, i.e.,  $\|c(\boldsymbol{z},t)\|\le C$ for  known  $C\ge 0$. Consider the RiSTL fragment
\begin{subequations}\label{eq:subclass}
\begin{align}
\psi \; &::= \; \top \; | \; \mu^{\text{Ch}} \; | \; \mu^{\text{Ri}}  \; | \; \psi' \wedge \psi'' \label{eq:psi_class}\\
\phi \; &::= \;  G_{[a,b]}\psi \; | \; F_{[a,b]} \psi \;|\; \psi'  \until{a}{b} \psi'' \; | \; \phi' \wedge \phi''\label{eq:phi_class}
\end{align}
\end{subequations}
where $\psi'$ and $\psi''$ are Boolean formulas of the form \eqref{eq:psi_class}, whereas $\phi'$ and $\phi''$ are of the form \eqref{eq:phi_class}. For specifications $\phi$ of the form \eqref{eq:subclass}, it holds that $(\boldsymbol{x},\boldsymbol{X},t)\models \phi$ if $\rho^\phi(\boldsymbol{x},\boldsymbol{X},t)\ge 0$, not requiring a strict inequality. The full RiSTL language as in \eqref{eq:full_RiSTL} can be dealt with when combining the proposed control laws  with timed automata theory  \cite{lindemann2019efficient}. Assume that the satisfaction of $\phi$ in \eqref{eq:subclass} depends on $\boldsymbol{x}$ and $\boldsymbol{X}$, but not on $\theta$. Assume also that $\phi$ consists of $M$ chance and risk predicates $\mu_m^{\text{Pr}}(\boldsymbol{x},\boldsymbol{X})$ and $\mu_m^{\text{Ri}}(\boldsymbol{x},\boldsymbol{X})$ for $m\in\{1,\hdots,M\}$ with associated predicate functions $h_m(\boldsymbol{x},\boldsymbol{X})$. 
\begin{assumption}
Each $h_m(\boldsymbol{x},\boldsymbol{X})$  is concave in $\boldsymbol{x}$.
\end{assumption} 
Let each $\mu_m^{\text{Pr}}(\boldsymbol{x},\boldsymbol{X})$ be associated with $\delta_m$ and each $\mu_m^{\text{Ri}}(\boldsymbol{x},\boldsymbol{X})$ be associated with $R_m(\cdot)$, $\beta_m$, and  $\gamma_m$.  We assume that  the mean $\boldsymbol{\tilde{\mu}}$, the covariance matrix $\tilde{\Sigma}$, and the probability density function $p_{\boldsymbol{X}}(\boldsymbol{X})$ of $\boldsymbol{X}$ is known.

\begin{problem}
Let $\phi$ be an RiSTL formula as in \eqref{eq:subclass}. Design a control law $\boldsymbol{u}(\boldsymbol{z},t)$ for  the system in \eqref{system_noise}  s.t. $\rho^\phi(\boldsymbol{x},\boldsymbol{X},t)\ge r\ge 0$, i.e., $(\boldsymbol{x},\boldsymbol{X},t)\models \phi$, where $r$ is maximized.
\end{problem}

\section{Proposed Problem Solution}
\label{sec:strategy}

Not that, for a fixed $\boldsymbol{x}$, each $h_m(\boldsymbol{x},\boldsymbol{X})$ has a mean  $\tilde{\mu}_{h_m}(\boldsymbol{x})$ and a variance $\tilde{\Sigma}_{h_m}(\boldsymbol{x})$. Let $p_{h_m}(h,\boldsymbol{x})$ denote the probability density function of $h_m(\boldsymbol{x},\boldsymbol{X})$ for a fixed $\boldsymbol{x}$.



\subsection{Determinization of RiSTL Specifications}
\label{sec:determ}

Note that $P(h_m(\boldsymbol{x},\boldsymbol{X})\ge 0)$ and $R(-h_m(\boldsymbol{x},\boldsymbol{X}))$ depend on $\boldsymbol{x}$. For given  $\delta_m,\beta_m\in(0,1)$ and $\gamma_m\in\mathbb{R}$,  define the sets
\begin{align*}
\mathfrak{X}_m^{\text{Ch}}(\delta_m)&:=\{\boldsymbol{x}\in\mathfrak{B}|P(h_m(\boldsymbol{x},\boldsymbol{X})\ge 0)\ge \delta_m\}\\
\mathfrak{X}_m^{\text{EV}}(\gamma_m)&:=\{\boldsymbol{x}\in\mathfrak{B}|E[-h_m(\boldsymbol{x},\boldsymbol{X})]\le \gamma_m\}\\
\mathfrak{X}_m^{\text{VaR}}(\beta_m,\gamma_m)&:=\{\boldsymbol{x}\in\mathfrak{B}|VaR_{\beta_m}(-h_m(\boldsymbol{x},\boldsymbol{X}))\le \gamma_m\}\\
\mathfrak{X}_m^{\text{CVaR}}(\beta_m,\gamma_m)&:=\{\boldsymbol{x}\in\mathfrak{B}|CVaR_{\beta_m}(-h_m(\boldsymbol{x},\boldsymbol{X}))\le \gamma_m\}
\end{align*}
where $\mathfrak{B}$ is an arbitrarily big compact and convex set, as further  explained in Section \ref{sec:control1}; $\mathfrak{X}_m^{\text{Ch}}(\delta_m)$ defines all $\boldsymbol{x}$ in $\mathfrak{B}$ for which the probability that $h_m(\boldsymbol{x},\boldsymbol{X})\ge 0$ is greater or equal than $\delta_m$, while $\mathfrak{X}_m^{\text{EV}}(\gamma_m)$, $\mathfrak{X}_m^{\text{VaR}}(\beta_m,\gamma_m)$, and $\mathfrak{X}_m^{\text{CVaR}}(\beta_m,\gamma_m)$ define all $\boldsymbol{x}$ in $\mathfrak{B}$ for which the EV, VaR, and CVaR of $-h_m(\boldsymbol{x},\boldsymbol{X})$  is less or equal than $\gamma_m$, respectively.  For a design parameter $c_m\in\mathbb{R}$, define
\begin{align*}
\mathfrak{X}_m(c_m):=\{\boldsymbol{x}\in\mathfrak{B}|h_m(\boldsymbol{x},\tilde{\boldsymbol{\mu}})-c_m\ge 0\}
\end{align*}
where the mean $\tilde{\boldsymbol{\mu}}$ of $\boldsymbol{X}$ is used. The set  $\mathfrak{X}_m(c_m)$ is compact and convex  since $h_m(\boldsymbol{x},\tilde{\boldsymbol{\mu}})$ is concave in $\boldsymbol{x}$.  If $\mathfrak{X}_m^{\text{Ch}}(\delta_m)\supseteq \mathfrak{X}_m(c_m)$, then $\boldsymbol{x}\in\mathfrak{X}_m(c_m)$  implies $\boldsymbol{x}\in{\mathfrak{X}}_m^{\text{Ch}}(\delta_m)$ (similarly for $\mathfrak{X}_m^{\text{EV}}(\gamma_m)$, $\mathfrak{X}_m^{\text{VaR}}(\beta_m,\gamma_m)$, and $\mathfrak{X}_m^{\text{CVaR}}(\beta_m,\gamma_m)$). In this case, an RiSTL formula can be determinized into an STL formula using $h_m(\boldsymbol{x},\tilde{\boldsymbol{\mu}})-c_m$ instead of $P(h_m(\boldsymbol{x},\boldsymbol{X})\ge 0)-\delta_m$ and $\gamma_m-R(-h_m(\boldsymbol{x},\boldsymbol{X}))$, mapping the stochastic control problem into a deterministic one (see Section \ref{sec:control1}).
\begin{example}\label{ex:1}
\begin{figure}
		\centering
		\includegraphics[scale=0.47]{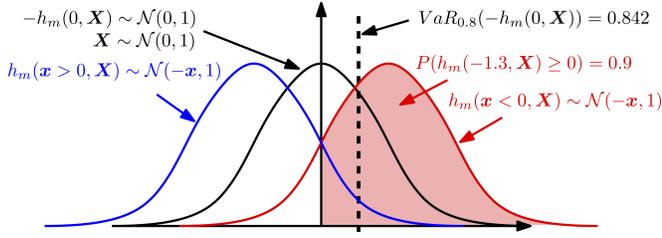}
		\caption{Illustrating Example \ref{ex:1} and the calculation of $c_m$.}
		\label{fig:exx2}
\end{figure}
Consider $\boldsymbol{X}\sim\mathcal{N}(0,1)$ and the predicate $\mu_m(\boldsymbol{x},\boldsymbol{X})$ with predicate function $h_m(\boldsymbol{x},\boldsymbol{X}):=\boldsymbol{X}-\boldsymbol{x}$. It holds that $\mathfrak{X}_m(c_m):=(-\infty,-c_m]$ since $\tilde{\boldsymbol{\mu}}=0$. Note that $h_m(\boldsymbol{x},\boldsymbol{X})\sim\mathcal{N}(-\boldsymbol{x},1)$ as in Fig. \ref{fig:exx2} where smaller $\boldsymbol{x}$ lead to larger $P(h_m(\boldsymbol{x},\boldsymbol{X})\ge 0)$ (indicated by the red area under the red curve). Consequently,  $\mathfrak{X}_m^{\text{Ch}}(\delta_m)=(-\infty,\boldsymbol{x}']$ where: 1) $\boldsymbol{x}'=0$ if $\delta_m=0.5$, 2) $\boldsymbol{x}'<0$ if $\delta_m>0.5$, and 3) $\boldsymbol{x}'>0$ if $\delta_m<0.5$. The idea is then, for a given $\delta_m$, to find $c_m$ such that $\mathfrak{X}_m(c_m)\subseteq\mathfrak{X}_m^{\text{Ch}}(\delta_m)$, e.g., for $\delta_m:=0.9$ it holds that $\boldsymbol{x}'=-1.3$ (see Fig. \ref{fig:exx2}) so that $c_m\ge 1.3$. Similarly, let  $\beta_m:=0.8$ so that $\mathfrak{X}_m^{\text{VaR}}(0.8,\gamma_m)=(-\infty,\boldsymbol{x}']$ where: 1) $\boldsymbol{x}'=0$ if $\gamma_m=0.842$, 2) $\boldsymbol{x}'<0$ if $\gamma_m<0.842$, and 3) $\boldsymbol{x}'>0$ if $\gamma_m>0.842$. For $\gamma_m=0.842$, it holds that $\boldsymbol{x}'=0$ so that $c_m\ge 0$ achieves $\mathfrak{X}_m(c_m)\subseteq\mathfrak{X}_m^{\text{VaR}}(\beta_m,\gamma_m)$.
\end{example}

Checking these set inclusions is, in general, nonconvex. When $h_m(\boldsymbol{x},\boldsymbol{X})$ is linear in $\boldsymbol{x}$, we can obtain the next result.
\begin{lemma}\label{lemma0}
Assume $h_m(\boldsymbol{x},\boldsymbol{X})=\boldsymbol{v}^T\boldsymbol{x}+h'(\boldsymbol{X})$ for  $\boldsymbol{v}\in\mathbb{R}^n$ and for  $h':\mathbb{R}^{\tilde{n}}\to\mathbb{R}$, then
\begin{align*}   
\mathfrak{X}_m^{\text{Ch}}(\delta_m)\supseteq \mathfrak{X}_m(c_m)  &\text{ iff } P(h_m(\boldsymbol{x}^*,\boldsymbol{X})\ge 0)\ge \delta_m\\
 \mathfrak{X}_m^{\text{EV}}(\gamma_m)\supseteq \mathfrak{X}_m(c_m)&\text{ iff } E[-h_m(\boldsymbol{x}^*,\boldsymbol{X})]\le \gamma_m \\
\mathfrak{X}_m^{\text{VaR}}(\beta_m,\gamma_m)\supseteq \mathfrak{X}_m(c_m) &\text{ iff } VaR_{\beta_m}(-h_m(\boldsymbol{x}^*,\boldsymbol{X}))\le \gamma_m \\
\mathfrak{X}_m^{\text{CVaR}}(\beta_m,\gamma_m)\supseteq \mathfrak{X}_m(c_m) &\text{ iff } CVaR_{\beta_m}(-h_m(\boldsymbol{x}^*,\boldsymbol{X}))\le \gamma_m 
\end{align*}
where $\boldsymbol{x}^*:=\underset{\boldsymbol{x}\in\mathfrak{X}_m(c_m)}{\text{argmin}}\; \boldsymbol{v}^T\boldsymbol{x}$  (a convex program).
\end{lemma}

Note in particular that, for $h_m(\boldsymbol{x},\boldsymbol{X})$ as in Lemma \ref{lemma0},   $VaR_{\beta_m}(-h_m(\boldsymbol{x}^*,\boldsymbol{X}))$ and $CVaR_{\beta_m}(-h_m(\boldsymbol{x}^*,\boldsymbol{X}))$ can be efficiently computed  \cite[Thm. 1]{rockafellar2000optimization}.

The RiSTL formula $\phi$ is now translated into an STL formula $\varphi$  by replacing chance and risk predicate in $\phi$ by 
\begin{align}\label{eq:STL_predicate}
\mu_m^{\text{STL}}(\boldsymbol{x}(t),\tilde{\boldsymbol{\mu}}):=\begin{cases}
\top & \text{if } h_m(\boldsymbol{x}(t),\tilde{\boldsymbol{\mu}})-c_m\ge 0\\
\bot &\text{otherwise. }
\end{cases}
\end{align}
The semantics of the STL formula $\varphi$ are, besides the evaluation of predicates in \eqref{eq:STL_predicate}, the same as for the RiSTL formula $\phi$ \cite{maler2004monitoring}. We also define quantitative semantics $\rho^\varphi(\boldsymbol{x},\tilde{\boldsymbol{\mu}},t)$ by letting $\rho^{\mu_m^{\text{STL}}}(\boldsymbol{x},\tilde{\boldsymbol{\mu}},t):=h_m(\boldsymbol{x}(t),\tilde{\boldsymbol{\mu}})-c_m$ and then following the recursive definition of RiSTL as introduced in Section \ref{sec:RiSTL} \cite{donze2}.  The following assumption is necessary  for $\mathfrak{X}_m(c_m)$ to be non empty and hence for $\varphi$ to be satisfiable.\footnote{More formally, necessity holds when the ``$>$" are replaced with ``$\ge$". }

\begin{assumption}\label{ass1}
	For each $m\in\{1,\hdots,M\}$, there exists $\boldsymbol{x}\in\mathbb{R}^n$ so that $h_m(\boldsymbol{x},\tilde{\boldsymbol{\mu}})-c_m>0$. Furthermore, for each $\psi'\wedge\psi''$ in $\varphi$ (recall \eqref{eq:psi_class}), there exists  $\boldsymbol{x}\in\mathbb{R}^n$ so that $\rho^{\psi'\wedge\psi''}(\boldsymbol{x},\tilde{\boldsymbol{\mu}},0)>0$.
\end{assumption}

Note that Assumption \ref{ass1} can efficiently be checked since $h_m(\boldsymbol{x},\tilde{\boldsymbol{\mu}})$ is concave in $\boldsymbol{x}$. It can  always be satisfied by a sufficiently small $c_m$ and  hence poses an upper bound on $c_m$. The next assumption is sufficient to ensure soundness in the sense that $(\boldsymbol{x},\boldsymbol{X},0)\models \varphi$ implies $(\boldsymbol{x},\boldsymbol{X},0)\models \phi$.

\begin{assumption}\label{ass2}
	For each $m\in\{1,\hdots,M\}$, $\mathfrak{X}_m^{\text{Pr}}(\delta_m)\supseteq \mathfrak{X}_m(c_m)$, $\mathfrak{X}_m^{\text{EV}}(\gamma_m)\supseteq \mathfrak{X}_m(c_m)$, $\mathfrak{X}_m^{\text{VaR}}(\beta_m,\gamma_m)\supseteq \mathfrak{X}_m(c_m)$, or $\mathfrak{X}_m^{\text{CVaR}}(\beta_m,\gamma_m)\supseteq \mathfrak{X}_m(c_m)$ depending on the predicate.
\end{assumption}

Increasing $c_m$ shrinks the set $\mathfrak{X}_m(c_m)$  so that Assumption~\ref{ass2} (verifiable by Lemma \ref{lemma0})  poses a lower bound on $c_m$.  
\begin{theorem}\label{thm:0} 
	Let Assumption \ref{ass2} hold. If $\boldsymbol{x}:\mathbb{R}_{\ge 0}\to\mathbb{R}^n$ is such that $(\boldsymbol{x},\boldsymbol{X},0)\models \varphi$, it follows that $(\boldsymbol{x},\boldsymbol{X},0)\models \phi$.

\end{theorem}

Finding a set of $c_m$ that satisfies Assumptions \ref{ass1} and \ref{ass2} may induce conservatism since the level sets of $\mathfrak{X}_m(c_m)$ may not be aligned with the level sets of $\mathfrak{X}_m^{\text{Ch}}(\delta_m)$, $\mathfrak{X}_m^{\text{EV}}(\gamma_m)$, $\mathfrak{X}_m^{\text{VaR}}(\beta_m,\gamma_m)$, and  $\mathfrak{X}_m^{\text{CVaR}}(\beta_m,\gamma_m)$. The next result shows when such conservatism can be avoided and alleviates finding a set of $c_m$ that satisfy Assumptions \ref{ass1} and \ref{ass2}.
\begin{lemma}\label{lem:cal}
	Assume that $h_m(\boldsymbol{x},\boldsymbol{X})=\boldsymbol{v}^T\boldsymbol{x}+h'(\boldsymbol{X})$ for  $\boldsymbol{v}\in\mathbb{R}^n$ and for  $h':\mathbb{R}^{\tilde{n}}\to\mathbb{R}$. Then there exists a design parameter $c_m$ so that $\mathfrak{X}_m(c_m)=\mathfrak{X}_m^{\text{Ch}}(\delta_m)$, $\mathfrak{X}_m(c_m)=\mathfrak{X}_m^{\text{EV}}(\gamma_m)$, $\mathfrak{X}_m(c_m)=\mathfrak{X}_m^{\text{VaR}}(\delta_m)$, or $\mathfrak{X}_m(c_m)=\mathfrak{X}_m^{\text{CVaR}}(\delta_m)$.
	
\end{lemma}

\subsection{Control Barrier Functions for Unicycle Dynamics}
\label{sec:control1}
Theorem \ref{thm:0} allows to map the stochastic control problem into a deterministic one.
The proposed control method is based on time-varying control barrier functions  where a function $\mathfrak{b}(\boldsymbol{x},\tilde{\boldsymbol{\mu}},t)$ encodes the STL formula $\varphi$ \cite{lindemann2018control}. Given that Assumptions \ref{ass1} and \ref{ass2} hold, we impose conditions on the function $\mathfrak{b}(\boldsymbol{x},\tilde{\boldsymbol{\mu}},t)$ as in  \cite[Steps A, B, and C]{lindemann2018control} that account for the STL semantics of $\varphi$; \cite{lindemann2019decentralized} presents a formally correct procedure to construct such $\mathfrak{b}(\boldsymbol{x},\tilde{\boldsymbol{\mu}},t)$. Define $\mathfrak{C}(\tilde{\boldsymbol{\mu}},t):=\{\boldsymbol{x}\in\mathbb{R}^2|\mathfrak{b}(\boldsymbol{x},\tilde{\boldsymbol{\mu}},t)\ge 0\}$ and note that $({\boldsymbol{x}},\tilde{\boldsymbol{\mu}},0)\models\varphi$ if $\boldsymbol{x}(t)\in\mathfrak{C}(\tilde{\boldsymbol{\mu}},t)$ for all $t\ge 0$ according to \cite{lindemann2018control}; $\mathfrak{C}(\tilde{\boldsymbol{\mu}},t)$ is ensured to be bounded so that we 
let $\mathfrak{D}$ be an open and bounded set with $\mathfrak{D}$ containing $\mathfrak{C}(\tilde{\boldsymbol{\mu}},t)$ for all $t\ge 0$. It is also ensured that $\boldsymbol{x}(0)\in\mathfrak{C}(\tilde{\boldsymbol{\mu}},0)$. In \cite{lindemann2018control} and \cite{lindemann2019decentralized}, the function  $\mathfrak{b}(\boldsymbol{x},\tilde{\boldsymbol{\mu}},t)$ is concave in the first argument  and piecewise continuous in the third argument with discontinuities at times $ \{s_0:=0,s_1,\hdots,s_q\}$ for some finite $q$.  
\begin{theorem}\label{thm:1}
Let the design parameters $c_m$ be so that  Assumptions \ref{ass1} and \ref{ass2} hold and let $\mathfrak{b}(\boldsymbol{x},\tilde{\boldsymbol{\mu}},t)$  be constructed for $\varphi$ as in \cite{lindemann2019decentralized}. If, for $\alpha>0$ and for all $(\boldsymbol{z},t)\in\mathfrak{D}\times(s_j,s_{j+1})$,  there exists a continuous control law $\boldsymbol{u}(\boldsymbol{z},t)$ such that
\begin{align}\label{eq:barrier_condition1}
\begin{split}
&{\frac{\partial \mathfrak{b}(\boldsymbol{x},\tilde{\boldsymbol{\mu}},t)}{\partial \boldsymbol{z}}}(f(\boldsymbol{z})+g(\boldsymbol{z})\boldsymbol{u}(\boldsymbol{z},t))+\frac{\partial \mathfrak{b}(\boldsymbol{x},\tilde{\boldsymbol{\mu}},t)}{\partial t} \\
&\hspace{2cm} \ge -\alpha\mathfrak{b}(\boldsymbol{x},\tilde{\boldsymbol{\mu}},t)+\Big\|{\frac{\partial \mathfrak{b}(\boldsymbol{x},\tilde{\boldsymbol{\mu}},t)}{\partial \boldsymbol{z}}}\Big\|C,
\end{split}
\end{align}
then $\rho^\phi(\boldsymbol{x},\boldsymbol{X},0)\ge r$ for some $r\ge 0$, i.e., $(\boldsymbol{x},\boldsymbol{X},0)\models \phi$.
\end{theorem}

For  unicycle-like dynamics in \eqref{system_noise}, the constraint in \eqref{eq:barrier_condition1} may not be feasible in case that
\begin{align*}
{\frac{\partial \mathfrak{b}(\boldsymbol{x},\tilde{\boldsymbol{\mu}},t)}{\partial \boldsymbol{z}}}g(\boldsymbol{z})={\frac{\partial \mathfrak{b}(\boldsymbol{x},\tilde{\boldsymbol{\mu}},t)}{\partial x_x}}\cos(\theta) + {\frac{\partial \mathfrak{b}(\boldsymbol{x},\tilde{\boldsymbol{\mu}},t)}{\partial x_y}}\sin(\theta)
\end{align*}
is equal to the zero vector. We use a near-identity diffeomorphism as in \cite{olfatinear} by the coordinate transformation 
\begin{align*}
\boldsymbol{p}:=\boldsymbol{x}+lR(\theta)\boldsymbol{e}_1
\end{align*} 
where $l>0$ is a design parameter and where $R(\theta):=\begin{bmatrix}
\cos(\theta) & -\sin(\theta)\\
\sin(\theta) & \cos(\theta)
\end{bmatrix}$ and $\boldsymbol{e}_1:=\begin{bmatrix} 1 & 0 \end{bmatrix}^T$. Note that $\dot{\boldsymbol{x}}=f_{\boldsymbol{x}}(\boldsymbol{z})+R(\theta)\boldsymbol{e}_1u_1+ c_{\boldsymbol{x}}(\boldsymbol{z},t)$ so that we can derive that 
\begin{align*}
\dot{\boldsymbol{p}}=f_{\boldsymbol{x}}(\boldsymbol{z})+g_{\boldsymbol{p}}(\boldsymbol{z}) \boldsymbol{u}+ c_{\boldsymbol{x}}(\boldsymbol{z},t)
\end{align*}
where $g_{\boldsymbol{p}}(\boldsymbol{z}):=\begin{bmatrix}
\cos(\theta) & -\sin(\theta)l\\
\sin(\theta) & \cos(\theta)l
\end{bmatrix}$ has full rank \cite[Lem. 1]{olfatinear}. Consider next the modified predicate
\begin{align}\label{eq:STL_predicate_m}
\hspace{-0.25cm}\bar{\mu}_m^{\text{STL}}(\boldsymbol{x}(t),\tilde{\boldsymbol{\mu}}):=\begin{cases}
\top & \text{if } h_m(\boldsymbol{x}(t),\tilde{\boldsymbol{\mu}})-c_m- \chi_m\ge 0\\
\bot &\text{otherwise }
\end{cases}
\end{align}
for $\chi_m>0$. The STL formula $\varphi$ is now transformed into the STL formula $\bar{\varphi}$ by replacing each predicate ${\mu}_m^{\text{STL}}(\boldsymbol{x}(t),\tilde{\boldsymbol{\mu}})$ in $\varphi$ by $\bar{\mu}_m^{\text{STL}}(\boldsymbol{x}(t),\tilde{\boldsymbol{\mu}})$. We then choose a sufficiently small $\chi_m$ for each $m\in\{1,\hdots,M\}$ so that Assumption \ref{ass1} holds for the modified predicate function $h_m(\boldsymbol{p},\tilde{\boldsymbol{\mu}})-c_m- \chi_m$ and $\bar{\varphi}$ and then construct $\mathfrak{b}(\boldsymbol{p},\tilde{\boldsymbol{\mu}},t)$ for $\bar{\varphi}$ as in \cite{lindemann2019decentralized}. We remark that we do not induce any conservatism and that choosing such $\chi_m$ is always possible. Note  that each $h_m(\boldsymbol{x},\tilde{\boldsymbol{\mu}})$ is locally Lipschitz continuous with Lipschitz constant $L_h^m$ on the domain $\mathfrak{D}$ so that $|h_m(\boldsymbol{x},\tilde{\boldsymbol{\mu}})-h_m(\boldsymbol{p},\tilde{\boldsymbol{\mu}})|\le L_h^m l$. We then select $l$  such that  $l\le \chi_m \mathbin{/} L_h^m$ for each $m\in\{1,\hdots,M\}$. Consequently, $h_m(\boldsymbol{p},\tilde{\boldsymbol{\mu}})-c_m-\chi_m\ge 0$ implies  $h_m(\boldsymbol{x},\tilde{\boldsymbol{\mu}})-c_m\ge0$.

\begin{theorem}\label{thm:2}
Let the design parameters $c_m$ be so that  Assumptions \ref{ass1} and \ref{ass2} hold, set $\chi_m$ and $l$  as instructed above, and let $\mathfrak{b}(\boldsymbol{p},\tilde{\boldsymbol{\mu}},t)$ be constructed for $\bar{\varphi}$ as in \cite{lindemann2019decentralized}. If $\alpha$ is as in \cite[Lemma 4]{lindemann2019decentralized}, then the control law $\boldsymbol{u}({\boldsymbol{z},t})$ given as 
\begin{subequations}\label{eq:barrier_condition3}
\begin{align}
& \boldsymbol{u}({\boldsymbol{z},t}):=\underset{\boldsymbol{u}\in\mathbb{R}^2}{\text{argmin}}\; \boldsymbol{u}^T\boldsymbol{u}\\
\begin{split}
&{\frac{\partial \mathfrak{b}(\boldsymbol{p},\tilde{\boldsymbol{\mu}},t)}{\partial \boldsymbol{p}}}(f_{\boldsymbol{x}}(\boldsymbol{z})+g_{\boldsymbol{p}}(\boldsymbol{z})\boldsymbol{u})+\frac{\partial \mathfrak{b}(\boldsymbol{p},\tilde{\boldsymbol{\mu}},t)}{\partial t} \\
&\hspace{2cm} \ge -\alpha\mathfrak{b}(\boldsymbol{p},\tilde{\boldsymbol{\mu}},t)+\Big\|{\frac{\partial \mathfrak{b}(\boldsymbol{p},\tilde{\boldsymbol{\mu}},t)}{\partial \boldsymbol{p}}}\Big\|C,\label{eq:barrier_single1}
\end{split}
\end{align}
\end{subequations}
ensures $\rho^\phi(\boldsymbol{x},\boldsymbol{X},0)\ge r$ for some $r\ge 0$, i.e., $({\boldsymbol{x}},\boldsymbol{X},0)\models\phi$.
\end{theorem}

In order to maximize $r$, one can find the set of $c_m$ that results in the largest $r$. This may, however, result in a tedious search. Another idea, possibly not obtaining the best $r$ but maximizing $r$ to some extent, is to obtain a set of $c_m$ so that Assumptions \ref{ass1} and \ref{ass2} hold and, instead of \eqref{eq:barrier_condition3}, solve
\begin{subequations}\label{eq:barrier_condition4}
\begin{align}
& \big(\boldsymbol{u}({\boldsymbol{z},t}),\epsilon({\boldsymbol{z},t})\big):=\underset{(\boldsymbol{u},\epsilon)\in\mathbb{R}^{2}\times\mathbb{R}_{\ge 0}}{\text{argmin}}\; \boldsymbol{u}^T\boldsymbol{u}-\epsilon\label{eq:rhd_robust1}\\
\begin{split}\label{eq:rhd_robust}
&{\frac{\partial \mathfrak{b}(\boldsymbol{p},\tilde{\boldsymbol{\mu}},t)}{\partial \boldsymbol{p}}}(f_{\boldsymbol{x}}(\boldsymbol{z})+g_{\boldsymbol{p}}(\boldsymbol{z})\boldsymbol{u})+\frac{\partial \mathfrak{b}(\boldsymbol{p},\tilde{\boldsymbol{\mu}},t)}{\partial t} \\
&\hspace{1.5cm} \ge -\alpha\mathfrak{b}(\boldsymbol{p},\tilde{\boldsymbol{\mu}},t)+\Big\|{\frac{\partial \mathfrak{b}(\boldsymbol{p},\tilde{\boldsymbol{\mu}},t)}{\partial \boldsymbol{p}}}\Big\|C+\epsilon.
\end{split}
\end{align}
\end{subequations}
Note that \eqref{eq:barrier_condition4}  is feasible for each $(\boldsymbol{z},t)\in\mathfrak{D}\times\mathbb{R}_{\ge 0}$ and, similarly to the proof of Theorem \ref{thm:2}, it can be shown that $\boldsymbol{u}({\boldsymbol{z},t})$ is continuous. Let $\epsilon_{\text{r}}:=\inf_{(\boldsymbol{z},t)\in\mathfrak{D}\times\mathbb{R}_{\ge 0}}\epsilon({\boldsymbol{z},t})$ and
\begin{align*}
\mathfrak{C}_{\text{r}}(\tilde{\boldsymbol{\mu}},t):=\{\boldsymbol{p}\in\mathbb{R}^2|\mathfrak{b}(\boldsymbol{p},\tilde{\boldsymbol{\mu}},t)\ge \epsilon_{\text{r}} \mathbin{/} \alpha\}.
\end{align*}
\begin{lemma}\label{thm3}
The control law $\boldsymbol{u}({\boldsymbol{z},t})$ in \eqref{eq:barrier_condition4} renders the set $\mathfrak{C}_{\text{r}}(\tilde{\boldsymbol{\mu}},t)$ forward invariant and attractive. 
\end{lemma}

Based on Lemma \ref{thm3}, we next find a lower bound on $r$. Therefore, we need the following result.

\begin{corollary}\label{corrrr}
The control law $\boldsymbol{u}(\boldsymbol{z},t)$ in \eqref{eq:barrier_condition4} results in $\rho^\varphi(\boldsymbol{x},\tilde{\boldsymbol{\mu}},0)\ge \epsilon_{\text{r}} \mathbin{/} \alpha$ if $\boldsymbol{p}(0)\in\mathfrak{C}_{\text{r}}(\tilde{\boldsymbol{\mu}},0)$.
\end{corollary}

Let us next define 
\begin{align*}
\mathfrak{X}_m^{\text{r}}(c_m):=\{\boldsymbol{x}\in\mathfrak{B}|h_m(\boldsymbol{x},\tilde{\boldsymbol{\mu}})-c_m\ge \epsilon_{\text{r}} \mathbin{/} \alpha\}
\end{align*}
for which $\mathfrak{X}_m^{\text{r}}(c_m)\subseteq \mathfrak{X}_m(c_m)$ with strict inclusion if $\epsilon_{\text{r}}>0$. Let	$r:=\min(r_1,\hdots, r_{M})$ with
$r_m:=\sup_{r\in\mathbb{R}} r \text{ s.t. CS}$ where $CS\in\{\mathfrak{X}_m^{\text{Ch}}(r+\delta_m)\supseteq \mathfrak{X}_m^{\text{r}}(c_m),
\mathfrak{X}_m^{\text{EV}}(\gamma_m-r)\supseteq \mathfrak{X}_m^{\text{r}}(c_m),
\mathfrak{X}_m^{\text{VaR}}(\beta_m,\gamma_m-r)\supseteq \mathfrak{X}_m^{\text{r}}(c_m),
\mathfrak{X}_m^{\text{CVaR}}(\beta_m,\gamma_m-r)\supseteq \mathfrak{X}_m^{\text{r}}(c_m)\}$ depending on the type of the predicate. The next result follows by  Corollary \ref{corrrr} and the definitions of $\mathfrak{X}_m^{\text{r}}(c_m)$, $r$, and the quantitative semantics of  $\phi$.

\begin{theorem}
The control law $\boldsymbol{u}(\boldsymbol{z},t)$ in \eqref{eq:barrier_condition4} results in $\rho^\phi(\boldsymbol{x},\tilde{\boldsymbol{\mu}},0)\ge r$.
\end{theorem}

\section{Simulations}
\label{sec:simulations}

\begin{figure}
        \centering
        \includegraphics[scale=0.33]{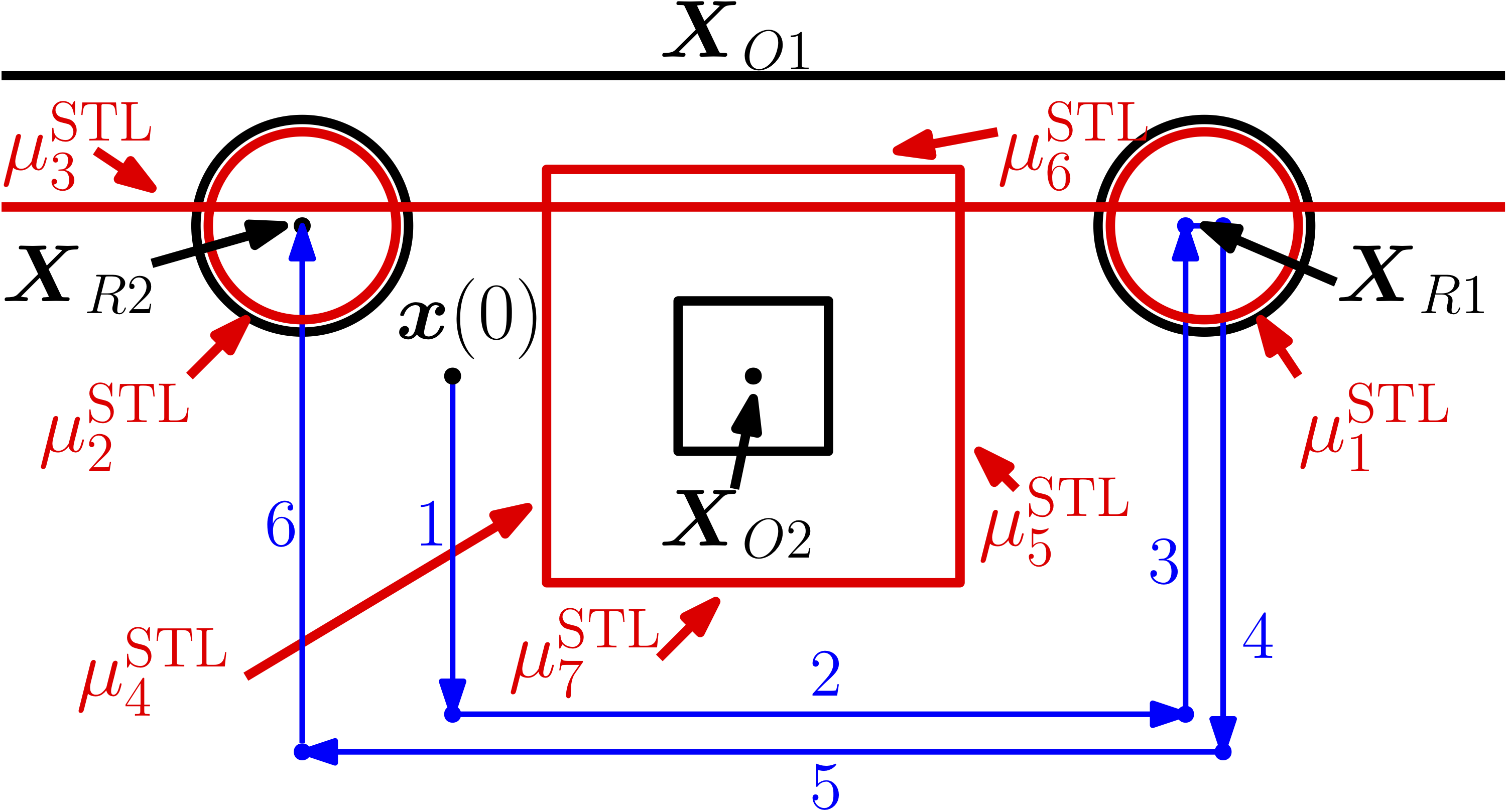}
        \caption{The RiSTL task $\phi$ contains the black chance and risk predicates $\mu_1^\text{Pr},\mu_2^\text{Pr},\mu_3^\text{Ri},\hdots,\mu_7^\text{Ri}$ . The determinized STL predicates $\mu_1^\text{STL},\hdots,\mu_7^\text{STL}$, contained in $\varphi$, are shown in red.}
        \label{fig:sim_overview}
\end{figure}
\begin{figure}
        \centering
        \includegraphics[scale=0.4]{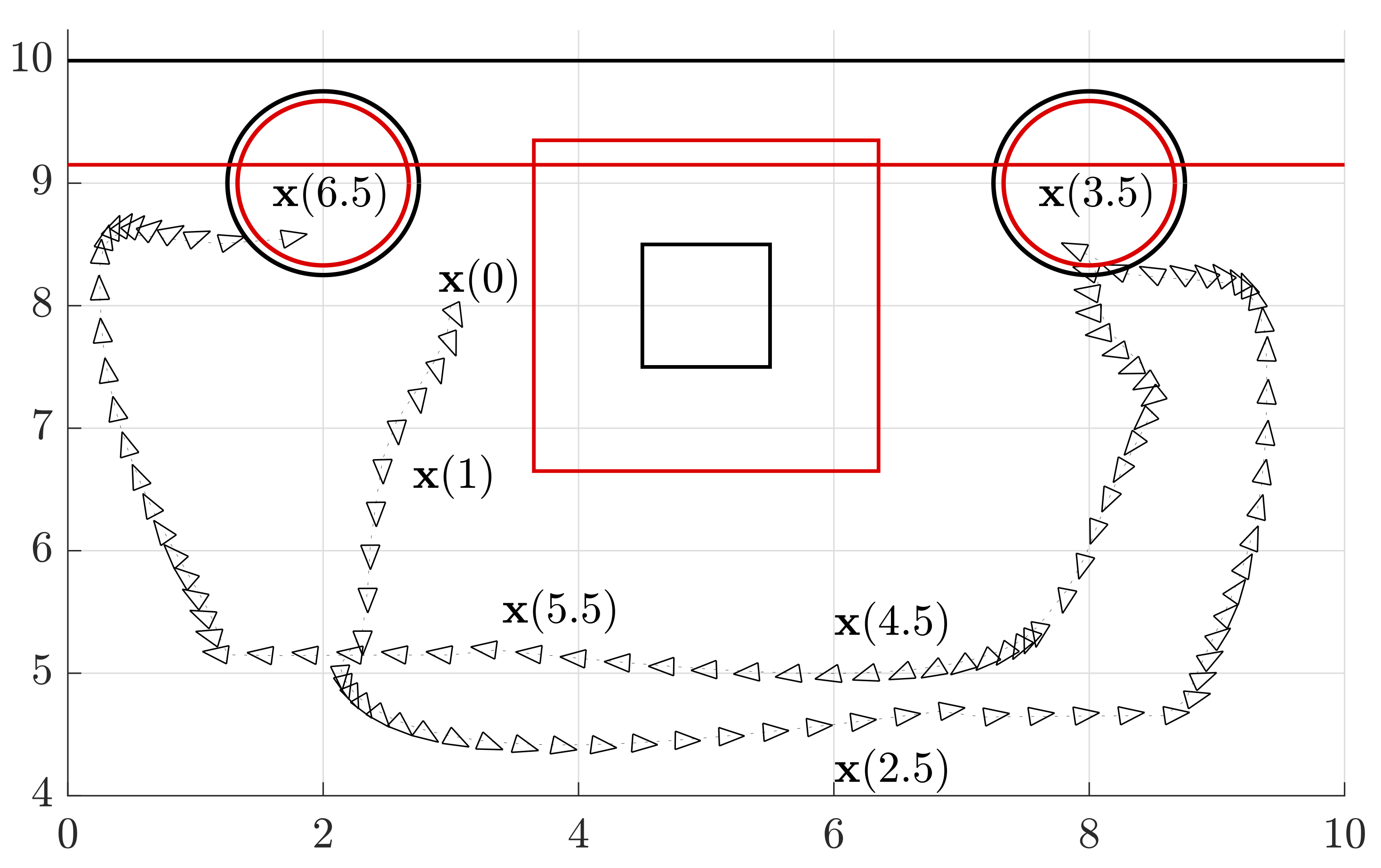}
        \caption{Simulation results.}
        \label{fig:sim}
\end{figure}

Consider the dynamics in \eqref{system_noise} with $f(\boldsymbol{z}):=\boldsymbol{0}$ and $c(\boldsymbol{z}):=0.5\cdot\begin{bmatrix} -\text{sat}(x_x) & -\text{sat}(x_y) & 0 \end{bmatrix}^T$ with $\boldsymbol{x}:=\begin{bmatrix}
x_x & x_y
\end{bmatrix}^T\in\mathbb{R}^2$ and where $\text{sat}(x)= x$ if $|x|\le 1$ and $\text{sat}(x)=1$ otherwise so that $C=0.5$. Furthemore, let ${\boldsymbol{X}}:=\begin{bmatrix}
{\boldsymbol{X}}_{O1}^T & {\boldsymbol{X}}_{O2}^T & {\boldsymbol{X}}_{R1}^T & {\boldsymbol{X}}_{R2}^T
\end{bmatrix}=\begin{bmatrix}
X_{O1,y} & X_{O2,x} & X_{O2,y}  & X_{R1,x} & X_{R1,y} & X_{R2,x} & X_{R2,y}
\end{bmatrix}^T$ where ${\boldsymbol{X}} \sim\mathcal{N}(\tilde{\boldsymbol{\mu}},\tilde{\Sigma})$ with $\tilde{\boldsymbol{\mu}}:=\begin{bmatrix}
10 & 5 & 8 & 8 & 9 & 2 & 9
\end{bmatrix}^T$ and $\tilde{\Sigma}:=\text{diag}(0.1,0.1,0.1,0.1,0.1,0.05,0.05)$. Let 
\begin{align*}
\phi:=  F_{[0,5]}(\phi_{R1}^\text{Pr} \wedge F_{[0,4]} \phi_{R2}^\text{Pr}) \wedge G_{[0,9]}  (\phi_{O1}^\text{Ri}\wedge  \phi_{O2}^\text{Ri})
\end{align*}
where $\phi_{R1}^\text{Pr}:=\mu_1^\text{Pr}$ and $\phi_{R1}^\text{Pr}:=\mu_2^\text{Pr}$ denote the probability of reaching regions indicated by ${\boldsymbol{X}}_{R1}$ and ${\boldsymbol{X}}_{R2}$ and $\phi_{O1}^\text{Ri}:=\mu_3^\text{Ri}$ and $\phi_{O2}^\text{Ri}:=\mu_4^\text{Ri}\vee\mu_5^\text{Ri}\vee\mu_6^\text{Ri}\vee \mu_7^\text{Ri}$ encode the risk of colliding with obstacles indicated by ${\boldsymbol{X}}_{O1}$ and ${\boldsymbol{X}}_{O2}$.\footnote{In particular, for $\epsilon_1:=0.75$ and $\epsilon_2:=0.5$ we define $h_1(\boldsymbol{x},{\boldsymbol{X}}):=\epsilon_1-\|\boldsymbol{x}-{\boldsymbol{X}}_{R1}\|$, $h_2(\boldsymbol{x},{\boldsymbol{X}}):=\epsilon_1-\|\boldsymbol{x}-{\boldsymbol{X}}_{R2}\|$,  $h_3(\boldsymbol{x},{\boldsymbol{X}}):=-x_y+ X_{O1,y}$, $h_4(\boldsymbol{x},{\boldsymbol{X}}):=-x_x+ X_{O2,x}-\epsilon_2$, $h_5(\boldsymbol{x},{\boldsymbol{X}}):=x_x- X_{O2,x}-\epsilon_2$, $h_6(\boldsymbol{x},{\boldsymbol{X}}):=x_y- X_{O2,y}-\epsilon_2$, and $h_7(\boldsymbol{x},{\boldsymbol{X}}):=-x_y+ X_{O2,y}-\epsilon_2$.} We use $\delta_1:=\delta_2:=0.85$ and $\gamma_m:=0$ and $\beta_m:=0.9$ for each $m\in\{3,\hdots,7\}$ while using CVaR. We obtain $\mathfrak{X}_m^{\text{Ch}}(\delta_m)\supseteq \mathfrak{X}_m(c_m)$ for $m\in\{1,2\}$ if $c_m:=0.08$ and $\mathfrak{X}_m^{\text{CVaR}}(\beta_m,\gamma_m)\supseteq \mathfrak{X}_m(c_m)$   for $m\in\{3,\hdots,7\}$ if $c_m:=0.85$, see Fig. \ref{fig:sim_overview}. Note that the system is never allowed to go around the obstacle $\boldsymbol{X}_{O2}$ from above which is deemed too risky. The RiSTL task $\phi$ and hence the STL tasks $\varphi$ and $\bar{\varphi}$ can not be encoded using the fragment in \eqref{eq:subclass}. We instead use the framework in \cite{lindemann2019efficient} to decompose $\bar{\varphi}$ into subtasks $\bar{\varphi}_i:=G_{[0,b_i]}\mu_{\text{inv},i}^{\text{STL}} \wedge F_{[b_i]} \mu_{\text{reach},i}^{\text{STL}}$ for $i\in\{1,\hdots,6\}$ with $b_i<b_{i+1}$  where each $\bar{\varphi}_i$ can be encoded in  \eqref{eq:subclass}. Sequentially satisfying each $\bar{\varphi}_i$ guarantees satisfaction of $\bar{\varphi}$, and consequently satisfaction of $\phi$. For an initial condition $\boldsymbol{x}(0)\models \bar{\mu}_3^{\text{STL}} \wedge \bar{\mu}_4^{\text{STL}} \wedge \neg \bar{\mu}_5^{\text{STL}} \wedge \neg \bar{\mu}_6^{\text{STL}} \wedge \neg \bar{\mu}_7^{\text{STL}}=:\mu_{\text{init}}^{\text{STL}}$, the algorithm in \cite{lindemann2019efficient} provides the sequence indicated by the blue waypoints in Fig. \ref{fig:sim_overview}. For instance, the first trajectory is constrained by $\mu_{\text{inv},1}^{\text{STL}}:=\mu_{\text{init}}^{\text{STL}}$,  $\mu_{\text{reach},1}^{\text{STL}}:=  \bar{\mu}_3^{\text{STL}} \wedge \bar{\mu}_4^{\text{STL}} \wedge \neg \bar{\mu}_5^{\text{STL}} \wedge \neg \bar{\mu}_6^{\text{STL}} \wedge \bar{\mu}_7^{\text{STL}}$, and $b_1:=1.5$. The simulation results are shown in Fig. \ref{fig:sim}. For $i=1$, we obtain $\epsilon_r=2.65$ and it is visible that the system always tries to maximize $r$ and hence the distance to the obstacles (see for instance the trajectory from $5.5$ to $6.5$ sec).

\section{Conclusion}
\label{sec:conclusion}

We presented risk signal temporal logic (RiSTL) by composing risk metrics with stochastic predicates to quantify the risk by which a predicate is not satisfied. We then considered unicycle-like dynamics in uncertain environments and showed that the stochastic control problem can be transformed into a deterministic one, which we solved by using time-varying control barrier functions. 
%

 \appendix
 \emph{Proof of Lemma \ref{lemma0}:}
Note that $-h_m(\boldsymbol{x},\boldsymbol{X})$ has mean $\tilde{\mu}_{-h_m}(\boldsymbol{x}):=-\boldsymbol{v}^T\boldsymbol{x}-E[h'(\boldsymbol{X})]$ and variance $\tilde{\Sigma}_{-h_m}$, which is a function of $h'(\boldsymbol{X})$ that does not depend on $\boldsymbol{x}$. Let $F_{-h_m}(h,\boldsymbol{x})$ be the cumulative distribution function of $-h_m(\boldsymbol{x},\boldsymbol{X})$ with $F_{-h_m}(h,\boldsymbol{x})=P(-\boldsymbol{v}^T\boldsymbol{x}-h'(\boldsymbol{X})\le h)=P(-h'(\boldsymbol{X})\le h+\boldsymbol{v}^T\boldsymbol{x})=F_{-h'}(h+\boldsymbol{v}^T\boldsymbol{x})$ where $F_{-h'}(h+\boldsymbol{v}^T\boldsymbol{x})$ is the cumulative distribution function of $-h'(\boldsymbol{X})$. Hence, $p_{-h_m}(h,\boldsymbol{x})=p_{-h'}(h+\boldsymbol{v}^T\boldsymbol{x})$, i.e., $p_{-h_m}(h,\boldsymbol{x})$ is of the same type for each $\boldsymbol{x}$ only shifted by $\boldsymbol{v}^T\boldsymbol{x}$.

Ch) Note that $\mathfrak{X}_m^{\text{Ch}}(\delta_m)\supseteq \mathfrak{X}_m(c_m)$ if and only if $\underset{\boldsymbol{x}\in\mathfrak{X}_m(c_m)}{\text{min}}\;P(h_m(\boldsymbol{x},\boldsymbol{X})\ge 0)\ge \delta_m$  by definition. Now
\begin{align*}
&\underset{\boldsymbol{x}\in\mathfrak{X}_m(c_m)}{\text{argmin}}P(h_m(\boldsymbol{x},\boldsymbol{X})\ge 0)=\underset{\boldsymbol{x}\in\mathfrak{X}_m(c_m)}{\text{argmin}}P(-h_m(\boldsymbol{x},\boldsymbol{X})\le 0)\\
&=\underset{\boldsymbol{x}\in\mathfrak{X}_m(c_m)}{\text{argmin}}P(-h'(\boldsymbol{X})\le \boldsymbol{v}^T\boldsymbol{x})=\underset{\boldsymbol{x}\in\mathfrak{X}_m(c_m)}{\text{argmin}}F_{-h'}(\boldsymbol{v}^T\boldsymbol{x})\\
&\stackrel{(a)}{=}\underset{\boldsymbol{x}\in\mathfrak{X}_m(c_m)}{\text{argmin}}\boldsymbol{v}^T\boldsymbol{x}=\boldsymbol{x}^*
\end{align*}
where (a) holds since $F_{-h'}(\boldsymbol{v}^T\boldsymbol{x})$ is nondecreasing. Hence, $\underset{\boldsymbol{x}\in\mathfrak{X}_m(c_m)}{\text{min}}P(h_m(\boldsymbol{x},\boldsymbol{X})\ge 0)=P(h_m(\boldsymbol{x}^*,\boldsymbol{X})\ge 0)$.

EV) Note  that $\mathfrak{X}_m^{\text{EV}}(\gamma_m)\supseteq \mathfrak{X}_m(c_m)$ if and only if $\underset{\boldsymbol{x}\in\mathfrak{X}_m(c_m)}{\text{max}}E[-h_m(\boldsymbol{x},\boldsymbol{X})]\le \gamma_m$ by definition. Now
\begin{align*}
&\underset{\boldsymbol{x}\in\mathfrak{X}_m(c_m)}{\text{argmax}} E[-h_m(\boldsymbol{x},\boldsymbol{X})]=\underset{\boldsymbol{x}\in\mathfrak{X}_m(c_m)}{\text{argmax}} \int_{-\infty}^\infty h p_{-h_m}(h,\boldsymbol{x}) \mathrm{d}h\\
&=\underset{\boldsymbol{x}\in\mathfrak{X}_m(c_m)}{\text{argmax}} \int_{-\infty}^\infty h p_{-h'}(h+\boldsymbol{v}^T\boldsymbol{x}) \mathrm{d}h\stackrel{(b)}{=}\underset{\boldsymbol{x}\in\mathfrak{X}_m(c_m)}{\text{argmin}}\boldsymbol{v}^T\boldsymbol{x}=\boldsymbol{x}^*
\end{align*}
where (b) holds since $\boldsymbol{x}^*$ maximizes $p_{-h'}(h+\boldsymbol{v}^T\boldsymbol{x})$ for each $h\in(-\infty,\infty)$, which maximizes the integral. Consequently, $\underset{\boldsymbol{x}\in\mathfrak{X}_m(c_m)}{\text{max}}E[-h_m(\boldsymbol{x},\boldsymbol{X})]=E[-h_m(\boldsymbol{x}^*,\boldsymbol{X})]$.

VaR) Note that $\mathfrak{X}_m^{\text{VaR}}(\beta_m,\gamma_m)\supseteq \mathfrak{X}_m(c_m)$ if and only if $\underset{\boldsymbol{x}\in\mathfrak{X}_m(c_m)}{\text{max}}VaR_{\beta_m}(-h_m(\boldsymbol{x},\boldsymbol{X}))\le \gamma_m$ by definition. Now
\begin{align*}
&\underset{\boldsymbol{x}\in\mathfrak{X}_m(c_m)}{\text{argmax}} VaR_{\beta_m}(-h_m(\boldsymbol{x},\boldsymbol{X}))\\
&=\underset{\boldsymbol{x}\in\mathfrak{X}_m(c_m)}{\text{argmax}}\min(d\in\mathbb{R}|P(-h_m(\boldsymbol{x},\boldsymbol{X})\le d)\ge \beta_m)\\
&=\underset{\boldsymbol{x}\in\mathfrak{X}_m(c_m)}{\text{argmax}}\min(d\in\mathbb{R}|F_{-h'}( d+\boldsymbol{v}^T\boldsymbol{x})\ge \beta_m)\\
&=\underset{\boldsymbol{x}\in\mathfrak{X}_m(c_m)}{\text{argmax}} h^*_0-\boldsymbol{v}^T\boldsymbol{x}=\underset{\boldsymbol{x}\in\mathfrak{X}_m(c_m)}{\text{argmin}}\boldsymbol{v}^T\boldsymbol{x}=\boldsymbol{x}^*
\end{align*}
where $h^*_0:=\min(d\in\mathbb{R}|F_{-h'}(d)\ge \beta_m)$. Hence, $\underset{\boldsymbol{x}\in\mathfrak{X}_m(c_m)}{\text{max}}VaR_{\beta_m}(-h_m(\boldsymbol{x},\boldsymbol{X}))=VaR_{\beta_m}(-h_m(\boldsymbol{x}^*,\boldsymbol{X}))$.

CVaR) Note that $\mathfrak{X}_m^{\text{CVaR}}(\beta_m,\gamma_m)\supseteq \mathfrak{X}_m(c_m)$ if and only if $\underset{\boldsymbol{x}\in\mathfrak{X}_m(c_m)}{\text{max}}CVaR_{\beta_m}(-h_m(\boldsymbol{x},\boldsymbol{X}))\le \gamma_m$. Noting $p_{-h_m}(h,\boldsymbol{x})=p_{-h'}(h+\boldsymbol{v}^T\boldsymbol{x})$, it holds that $\underset{\boldsymbol{x}\in\mathfrak{X}_m(c_m)}{\text{argmax}}CVaR_{\beta_m}(-h_m(\boldsymbol{x},\boldsymbol{X}))=\underset{\boldsymbol{x}\in\mathfrak{X}_m(c_m)}{\text{argmax}}VaR_{\beta_m}(-h_m(\boldsymbol{x},\boldsymbol{X}))=\boldsymbol{x}^*$. Consequently, $\underset{\boldsymbol{x}\in\mathfrak{X}_m(c_m)}{\text{max}}CVaR_{\beta_m}(-h_m(\boldsymbol{x},\boldsymbol{X}))=CVaR_{\beta_m}(-h_m(\boldsymbol{x}^*,\boldsymbol{X}))$.

\emph{Proof of Theorem \ref{thm:0}:}		Due to Assumption \ref{ass2}, $\boldsymbol{x}\in\mathfrak{X}_m(c_m)$ implies $\boldsymbol{x}\in\mathfrak{X}_m^{\text{Pr}}(\delta_m)$, $\boldsymbol{x}\in\mathfrak{X}_m^{\text{EV}}(\gamma_m)$, $\boldsymbol{x}\in\mathfrak{X}_m^{\text{VaR}}(\beta_m,\gamma_m)$, or $\boldsymbol{x}\in\mathfrak{X}_m^{\text{CVaR}}(\beta_m,\gamma_m)$. Since the semantics of STL and RiSTL only differ on the predicate level and $\phi$ does not contain negations,  $(\boldsymbol{x},\boldsymbol{X},0)\models \varphi$ implies $(\boldsymbol{x},\boldsymbol{X},0)\models \phi$.

\emph{Proof of Lemma \ref{lem:cal}:}			The proof follows by noting that all $\boldsymbol{x}$ that satisfy $\boldsymbol{v}^T\boldsymbol{x}=\nu$ for some $\nu\in\mathbb{R}$ result in the same $P(h_m(\boldsymbol{x},\boldsymbol{X})\ge 0)$, $E[-h_m(\boldsymbol{x},\boldsymbol{X})]$, $VaR_\beta(-h_m(\boldsymbol{x},\boldsymbol{X}))$, and $CVaR_\beta(-h_m(\boldsymbol{x},\boldsymbol{X}))$, respectively. In other words, the level sets of $P(h_m(\boldsymbol{x},\boldsymbol{X})\ge 0)$, $E[-h_m(\boldsymbol{x},\boldsymbol{X})]$, $VaR_\beta(-h_m(\boldsymbol{x},\boldsymbol{X}))$, and $CVaR_\beta(-h_m(\boldsymbol{x},\boldsymbol{X}))$ form again hyperplanes with normal vector $\boldsymbol{v}$. Noting that the level sets of $\mathfrak{X}_m(c_m)$ also result in a hyperplane with normal vector $\boldsymbol{v}$ that can be shifted by $c_m$ completes the proof.

\emph{Proof of Theorem \ref{thm:1}:}
Recall that  $\boldsymbol{z}(t):=\begin{bmatrix} \boldsymbol{x}(t)^T & \theta(t)\end{bmatrix}^T$. Since $\boldsymbol{u}(\boldsymbol{z},t)$ is continuous, there exist solutions $\boldsymbol{z}:[0,\tau_{\text{max}})\to\mathfrak{D}$  to \eqref{system_noise} with $\tau_{\text{max}}>0$. Now, \eqref{eq:barrier_condition1} implies ${\frac{\partial \mathfrak{b}(\boldsymbol{x},\tilde{\boldsymbol{\mu}},t)}{\partial \boldsymbol{z}}}(f(\boldsymbol{z})+g(\boldsymbol{z})\boldsymbol{u}(\boldsymbol{z},t)+c(\boldsymbol{z},t))+\frac{\partial \mathfrak{b}(\boldsymbol{x},\tilde{\boldsymbol{\mu}},t)}{\partial t} \ge -\alpha\mathfrak{b}(\boldsymbol{x},\tilde{\boldsymbol{\mu}},t)$ so that, for all $t\in(0,\min(\tau_{\text{max}},s_1))$, $\dot{\mathfrak{b}}(\boldsymbol{x}(t),\tilde{\boldsymbol{\mu}},t)\ge-\alpha\mathfrak{b}(\boldsymbol{x}(t),\tilde{\boldsymbol{\mu}},t)$. Due to  \cite[Lem.~4.4]{Kha96}, the Comparison Lemma \cite[Ch.~3.4]{Kha96}, and since $\mathfrak{b}(\boldsymbol{x}(0),\tilde{\boldsymbol{\mu}},0)\ge 0$, it follows that $\mathfrak{b}(\boldsymbol{x}(t),\tilde{\boldsymbol{\mu}},t)\ge  0$, i.e.,  $\boldsymbol{x}(t)\in\mathfrak{C}(\tilde{\boldsymbol{\mu}},t)$, for all $t\in[0,\min(\tau_{\text{max}},s_1))$. If $\tau_{\text{max}}\ge s_1$, it holds $\boldsymbol{x}(t)\in\mathfrak{C}(\tilde{\boldsymbol{\mu}},t)$ for all $t\in[s_1,\min(\tau_{\text{max}},s_2))$. By \cite{lindemann2019decentralized}, for each $s_j$ with $j\in\{1,\hdots,q\}$, it holds that $\lim_{\tau\to s_j^-} \mathfrak{C}(\tilde{\boldsymbol{\mu}},\tau)\subseteq \mathfrak{C}(\tilde{\boldsymbol{\mu}},s_j)$ so that $\boldsymbol{x}(s_1)\in\mathfrak{C}(\tilde{\boldsymbol{\mu}},s_1)$. This argument can be repeated unless $\tau_{\text{max}}<s_j$ for some $j$; however, $\mathfrak{b}(\boldsymbol{x}(t),\tilde{\boldsymbol{\mu}},t)\ge 0$ implies that $\boldsymbol{x}(t)\in\mathfrak{B}$ for the compact set  $\mathfrak{B}\subset \mathbb{R}^n$ and for all $t\in[0,\tau_{\text{max}})$; $\theta(t)$ will evolve in a compact set since $f_\theta(\boldsymbol{z})$ and $c_\theta(\boldsymbol{z},t)$ are bounded so that  $\tau_{\text{max}}=\infty$ \cite[Thm.~3.3]{Kha96}. By \cite{lindemann2018control},  $\boldsymbol{x}(t)\in\mathfrak{C}(\tilde{\boldsymbol{\mu}},t)$ for all $t\ge 0$ so that $({\boldsymbol{x}},\tilde{\boldsymbol{\mu}},0)\models\varphi$, i.e., $\rho^\varphi(\boldsymbol{x},\tilde{\boldsymbol{\mu}},0)\ge r'$ for some $r'\ge 0$. Hence, $({\boldsymbol{x}},\tilde{\boldsymbol{\mu}},0)\models\phi$, i.e., $\rho^\phi(\boldsymbol{x},\tilde{\boldsymbol{\mu}},0)\ge r$ for some $r\ge 0$,  by Theorem \ref{thm:0}.

\emph{Proof of Theorem \ref{thm:2}:}  If $(\boldsymbol{z},t)\in\mathbb{R}^{3}\times(s_j,s_{j+1})$ with $\frac{\partial \mathfrak{b}(\boldsymbol{p},\tilde{\boldsymbol{\mu}},t)}{\partial \boldsymbol{p}}g_{\boldsymbol{p}}(\boldsymbol{z})\neq \boldsymbol{0}$, \eqref{eq:barrier_condition3} is feasible and $\boldsymbol{u}(\boldsymbol{z},t)$ is locally Lipschitz continuous at $(\boldsymbol{z},t)$ \cite[Thm. 8]{xu2015robustness}. Note that $\frac{\partial \mathfrak{b}(\boldsymbol{p},\tilde{\boldsymbol{\mu}},t)}{\partial \boldsymbol{p}}g_{\boldsymbol{p}}(\boldsymbol{z})=\boldsymbol{0}$ if and only if $\frac{\partial \mathfrak{b}(\boldsymbol{p},\tilde{\boldsymbol{\mu}},t)}{\partial \boldsymbol{p}}=\boldsymbol{0}$ since $g_{\boldsymbol{p}}(\boldsymbol{z})$ has full rank. If $(\boldsymbol{z},t)\in\mathbb{R}^{3}\times(s_j,s_{j+1})$ with $\frac{\partial \mathfrak{b}(\boldsymbol{p},\tilde{\boldsymbol{\mu}},t)}{\partial \boldsymbol{p}}= \boldsymbol{0}$, \eqref{eq:barrier_single1} is satisfied since $\frac{\partial \mathfrak{b}(\boldsymbol{p},\tilde{\boldsymbol{\mu}},t)}{\partial t}\ge -\alpha\mathfrak{b}(\boldsymbol{p},\tilde{\boldsymbol{\mu}},t)+\chi$ for some $\chi>0$  due the choice of $\alpha$ so that $\boldsymbol{u}(\boldsymbol{z},t):=\boldsymbol{0}$.\footnote{See \cite[Lemma 4]{lindemann2019decentralized} for details on the choice of $\alpha$. The intution is that $\mathfrak{b}(\boldsymbol{p},\tilde{\boldsymbol{\mu}},t)$ is concave in $\boldsymbol{p}$ so that $\frac{\partial \mathfrak{b}(\boldsymbol{p},\tilde{\boldsymbol{\mu}},t)}{\partial \boldsymbol{p}}= \boldsymbol{0}$ only happens at the global optimum $\boldsymbol{p}^*$ where $\mathfrak{b}(\boldsymbol{p}^*,\tilde{\boldsymbol{\mu}},t)>0$. Then choosing $\alpha$ large enough ensures that $\frac{\partial \mathfrak{b}(\boldsymbol{p}^*,\tilde{\boldsymbol{\mu}},t)}{\partial t}\ge -\alpha\mathfrak{b}(\boldsymbol{p}^*,\tilde{\boldsymbol{\mu}},t)+\chi$ for some $\chi>0$.} Due to continuity of $\frac{\partial \mathfrak{b}(\boldsymbol{p},\tilde{\boldsymbol{\mu}},t)}{\partial t}$ and $\alpha\mathfrak{b}(\boldsymbol{p},\tilde{\boldsymbol{\mu}},t)$, there exists a neighborhood $\mathcal{U}$ around $(\boldsymbol{p},t)$  so that, for each $(\boldsymbol{p}',t')\in\mathcal{U}$, $\frac{\partial \mathfrak{b}(\boldsymbol{p}',\tilde{\boldsymbol{\mu}},t')}{\partial t}\ge -\alpha\mathfrak{b}(\boldsymbol{p}',\tilde{\boldsymbol{\mu}},t'))$ and consequently $\boldsymbol{u}(\boldsymbol{p}',t')=\boldsymbol{0}$. Hence, $\boldsymbol{u}(\boldsymbol{z},t)$  is continuous on $\mathbb{R}^{3}\times(s_j,s_{j+1})$ so that, similarly to the proof of Theorem~\ref{thm:1},   $\rho^{\bar{\varphi}}(\boldsymbol{p},\tilde{\boldsymbol{\mu}},0) \ge 0$ which implies $\rho^{\varphi}(\boldsymbol{p},\tilde{\boldsymbol{\mu}},0) \ge \min(\chi_1,\hdots,\chi_{M})$ by \eqref{eq:STL_predicate} and \eqref{eq:STL_predicate_m} and the syntax of $\phi$ (and consequently $\varphi$) in \eqref{eq:subclass} that exclude disjunctions and negations. By the choice of $l$, this implies $({\boldsymbol{x}},\tilde{\boldsymbol{\mu}},0)\models{\varphi}$ so that again $({\boldsymbol{x}},\boldsymbol{X},0)\models\phi$, i.e., $\rho^\phi(\boldsymbol{x},\boldsymbol{X},0)\ge r$ for some $r\ge 0$, as in proof of Theorem \ref{thm3}.

\emph{Proof of Lemma \ref{thm3}:}
First note that, for each solution $\boldsymbol{p}:\mathbb{R}_{\ge 0}\to \mathbb{R}^n$ that arises under $\boldsymbol{u}({\boldsymbol{z},t})$,  $\dot{\mathfrak{b}}(\boldsymbol{p}(t),\tilde{\boldsymbol{\mu}},t)\ge -\alpha\mathfrak{b}(\boldsymbol{p}(t),\tilde{\boldsymbol{\mu}},t)+\epsilon_{\text{r}}$ due to \eqref{eq:rhd_robust}. For $\epsilon_{\text{r}}\ge 0$ and $\alpha> 0$, the initial value problem $\dot{v}(t)=-\alpha v(t)+\epsilon_{\text{r}}$ with $v(0)\ge \epsilon_{\text{r}} \mathbin{/} \alpha$ has the solution  $v(t)= \exp(-\alpha t) (v(0)-\epsilon_{\text{r}} \mathbin{/} \alpha)+\epsilon_{\text{r}} \mathbin{/} \alpha\ge \epsilon_{\text{r}} \mathbin{/} \alpha$. By the Comparison Lemma \cite[Ch.~3.4]{Kha96}, it follows that ${\mathfrak{b}}(\boldsymbol{p}(t),\tilde{\boldsymbol{\mu}},t)\ge v(t)\ge \epsilon_{\text{r}} \mathbin{/} \alpha$ so that $\boldsymbol{p}(t)\in\mathfrak{C}_{\text{r}}(\tilde{\boldsymbol{\mu}},t)$ for all $t\ge 0$ if $\boldsymbol{p}(0)\in\mathfrak{C}_{\text{r}}(\tilde{\boldsymbol{\mu}},0)$.  If $v(0)\in\mathbb{R}$, $v(t)\to \epsilon \mathbin{/} \alpha$ as $t\to\infty$ so that  attractivity of $\mathfrak{C}_{\text{r}}(\tilde{\boldsymbol{\mu}},t)$ follows under $\boldsymbol{u}({\boldsymbol{z},t})$ by again using the Comparison Lemma. 

\emph{Proof of Corollary \ref{corrrr}:}
Due to Lemma \ref{thm3} and since $\boldsymbol{p}(0)\in\mathfrak{C}_{\text{r}}(\tilde{\boldsymbol{\mu}},0)$,  it holds that $\mathfrak{b}(\boldsymbol{p}(t),\tilde{\boldsymbol{\mu}},t)\ge \epsilon_{\text{r}} \mathbin{/} \alpha$ for all $t\ge 0$. By the construction of $\mathfrak{b}(\boldsymbol{p},\tilde{\boldsymbol{\mu}},t)$ and \cite{lindemann2018control} it follows that $\rho^{\bar{\varphi}}(\boldsymbol{p},\tilde{\boldsymbol{\mu}},0)\ge \epsilon_{\text{r}} \mathbin{/} \alpha$ so that $\rho^{{\varphi}}(\boldsymbol{p},\tilde{\boldsymbol{\mu}},0)\ge \min(\chi_1,\hdots,\chi_{M})+\epsilon_{\text{r}} \mathbin{/} \alpha$. Since, for each $m\in\{1,\hdots,M\}$, $h_m(\boldsymbol{p},\tilde{\boldsymbol{\mu}})-c_m-\chi_m\ge  \epsilon_{\text{r}} \mathbin{/} \alpha$ implies  $h_m(\boldsymbol{x},\tilde{\boldsymbol{\mu}})-c_m> \epsilon_{\text{r}} \mathbin{/} \alpha$  by the choice of $l$, it follows that $\rho^{{\varphi}}(\boldsymbol{x},\tilde{\boldsymbol{\mu}},0)\ge \epsilon_{\text{r}} \mathbin{/} \alpha$.

\bibliographystyle{IEEEtran}
\bibliography{literature}

\newpage

\addtolength{\textheight}{-12cm}   

\end{document}